\documentclass[12pt,twoside]{amsart}
\usepackage[margin=3cm]{geometry}
\usepackage[colorlinks=false]{hyperref}
\usepackage[english]{babel}
\usepackage{graphicx,titling}
\usepackage{float}
\usepackage{amsmath,amsfonts,amssymb,amsthm}
\usepackage{lipsum}
\usepackage[T1]{fontenc}
\usepackage{fourier}
\usepackage{color}
\usepackage[latin1]{inputenc}
\usepackage{esint}
\usepackage{caption}
\usepackage{ccicons}

\makeatletter
\def\blfootnote{\xdef\@thefnmark{}\@footnotetext}
\makeatother

\newcommand\ccnote{
    \blfootnote{\copyright\,\, Quoc-Hung Nguyen, Matthew Rosenzweig, and Sylvia Serfaty}
    \blfootnote{\ccLogo\, \ccAttribution\,\, Licensed under a \href{https://creativecommons.org/licenses/by/4.0/}{Creative Commons Attribution License (CC-BY)}.}
}

\usepackage[export]{adjustbox}
\numberwithin{equation}{section}
\usepackage{setspace}
\renewcommand{\le}{\leqslant}
\renewcommand{\leq}{\leqslant}
\renewcommand{\ge}{\geqslant}
\renewcommand{\geq}{\geqslant}
\renewcommand{\mathbb}{\varmathbb}
\usepackage{fancyhdr}
\newtheorem{theorem}{Theorem}[section]

\newtheorem{corollary}[theorem]{Corollary}
\newtheorem{proposition}[theorem]{Proposition}
\newtheorem{remark}[theorem]{Remark}
\usepackage{listings}
\usepackage{lstautogobble}
\usepackage{enumerate}
\usepackage{thmtools}
\usepackage{thm-restate}
\usepackage{verbatim}

\usepackage{mathtools}
\usepackage{physics}
\usepackage{color}
\usepackage[capitalise]{cleveref}
\usepackage[bottom]{footmisc}
\crefformat{equation}{(#2#1#3)}

\usepackage{float}
\restylefloat{table}

\usepackage{mathrsfs}

\DeclarePairedDelimiter\ipp{\langle}{\rangle}

\DeclarePairedDelimiter{\paren}{\lparen}{\rparen}

\DeclareMathOperator{\supp}{supp}

\newcommand{\M}{{\mathcal{M}}}

\newcommand{\p}{{\partial}}

\renewcommand{\d}{\delta}

\newcommand{\R}{{\mathbb{R}}}

\newcommand{\N}{{\mathbb{N}}}

\newcommand{\K}{{\mathsf{K}}}

\renewcommand{\H}{{\mathcal{H}}}
\newcommand{\T}{{\mathbb{T}}}
\newcommand{\g}{{\mathsf{g}}}
\newcommand{\G}{{\mathsf{G}}}

\newcommand{\Sc}{{\mathcal{S}}}

\renewcommand{\M}{{\mathbb{M}}}

\renewcommand{\k}{\mathsf{k}}
\newcommand{\I}{\mathbb{I}}

\newcommand{\wt}{\widetilde}
\newcommand{\tl}{\tilde}

\newcommand{\D}{\Delta}

\newcommand{\ph}{\phantom{=}}
\newcommand{\nn}{\nonumber}

\newcommand{\ul}{\underline}

\newcommand{\ux}{\underline{x}}

\newcommand{\ep}{\varepsilon}
\newcommand{\al}{\alpha}
\newcommand{\be}{\beta}

\newcommand{\Dm}{|\nabla|}

\newcommand{\Fr}{{F}}
\renewcommand{\P}{\mathcal{P}}
\newcommand{\Te}{\mathrm{Term}}

\let\div\relax
\DeclareMathOperator{\div}{div}

\def\Xint#1{\mathchoice
{\XXint\displaystyle\textstyle{#1}}%
{\XXint\textstyle\scriptstyle{#1}}%
{\XXint\scriptstyle\scriptscriptstyle{#1}}%
{\XXint\scriptscriptstyle\scriptscriptstyle{#1}}%
\!\int}
\def\XXint#1#2#3{{\setbox0=\hbox{$#1{#2#3}{\int}$ }
\vcenter{\hbox{$#2#3$ }}\kern-.6\wd0}}

\def\dashint{\Xint-}

\def \be{\begin{equation}}
\def \ee{\end{equation}}
\def \hal{\frac{1}{2}}
\def\({\left(}
\def\){\right)}
\def \ep{\varepsilon}
\def\nab{\nabla}
\def\indic{\mathbf{1}}

\address{Quoc-Hung Nguyen, Academy of Mathematics and Systems Science, Chinese Academy of Sciences (CAS), Beijing, China}
\email{qhnguyen@amss.ac.cn}
\medskip
\address{Matthew Rosenzweig, Massachusetts Institute of Technology, Department of Mathematics, Cambridge, MA} 
\email{mrosenzw@mit.edu}
\medskip
\address{Sylvia Serfaty, Courant Institute of Mathematical Sciences, New York University, New York City, NY}
\email{serfaty@cims.nyu.edu}

\begin{document}

\thispagestyle{empty}

\begin{minipage}{0.28\textwidth}
\begin{figure}[H]
\includegraphics[width=2.5cm,height=2.5cm,left]{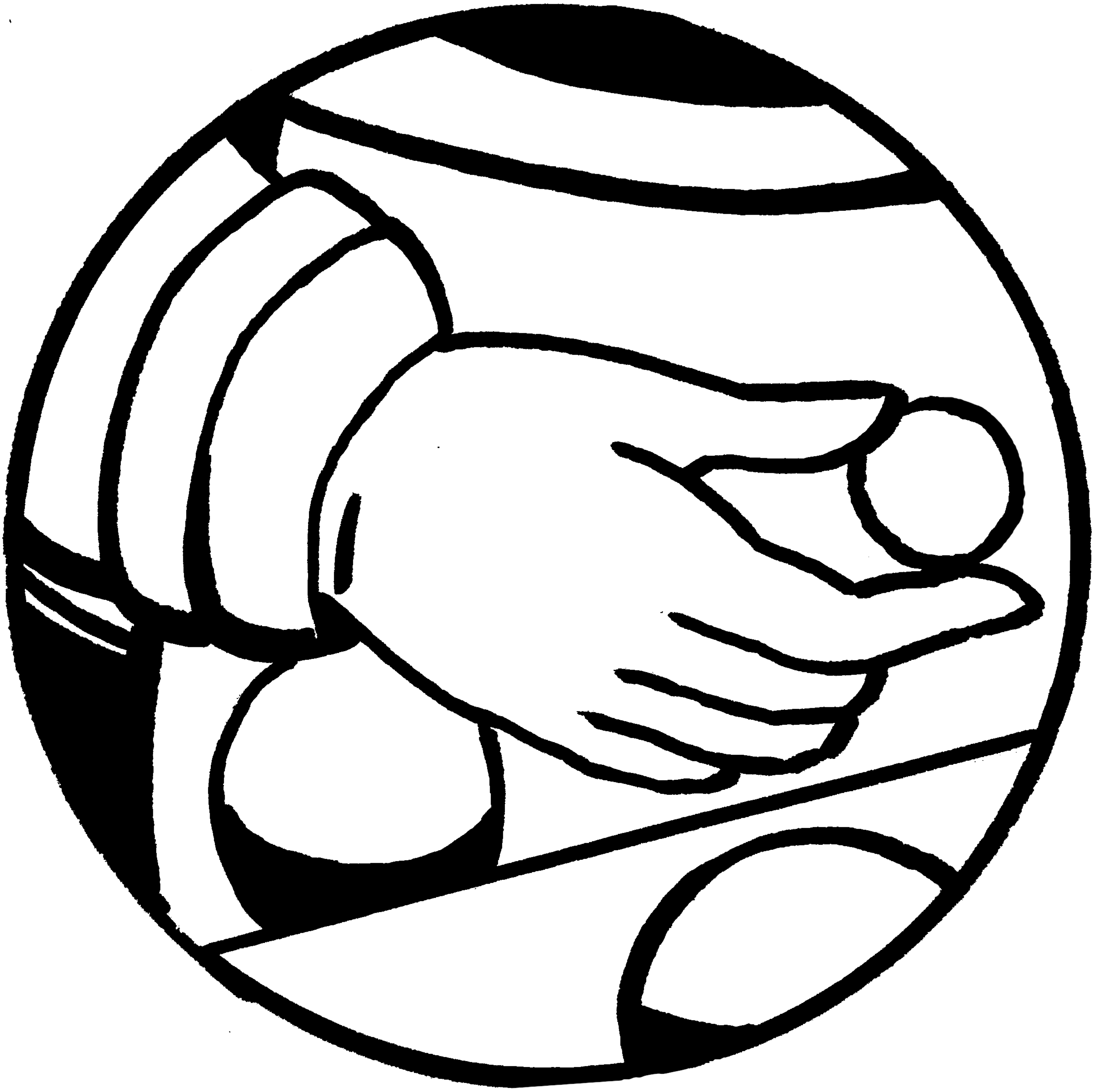}
\end{figure}
\end{minipage}
\begin{minipage}{0.7\textwidth} 
\begin{flushright}
Ars Inveniendi Analytica (2022), Paper No. 4, 45 pp.
\\
DOI 10.15781/nvv7-jy87
\\
ISSN: 2769-8505
\end{flushright}
\end{minipage}

\ccnote

\vspace{0.5cm}


\begin{center}
\begin{huge}
\textit{Mean-field limits \\ of Riesz-type singular flows}


\end{huge}
\end{center}

\vspace{1cm}


\begin{minipage}[t]{.28\textwidth}
\begin{center}
{\large{\bf{Quoc-Hung Nguyen}}} \\
\vskip0.15cm
\footnotesize{Academy of Mathematics and Systems Science, Chinese Academy of Sciences}
\end{center}
\end{minipage}
\hfill
\noindent
\begin{minipage}[t]{.28\textwidth}
\begin{center}
{\large{\bf{Matthew Rosenzweig}}} \\
\vskip0.15cm
\footnotesize{Massachusetts Institute of Technology}
\end{center}
\end{minipage}
\hfill
\noindent
\begin{minipage}[t]{.28\textwidth}
\begin{center}
{\large{\bf{Sylvia Serfaty}}} \\
\vskip0.15cm
\footnotesize{Courant Institute of Mathematical Sciences, New York University} 
\end{center}
\end{minipage}

\vspace{0.7cm}


\begin{center}
\noindent \em{Communicated by Jacob Bedrossian}
\end{center}

\vspace{0.7cm}


\noindent \textbf{Abstract.} \textit{We provide a  proof of  mean-field convergence of first-order dissipative or conservative dynamics  of particles with Riesz-type singular interaction (the model interaction is an inverse power $s$ of the distance for any $0<s<d$) when assuming a certain regularity of the solutions to the limiting evolution equations. It relies on a modulated-energy approach, as introduced in previous works where it was restricted to the Coulomb and super-Coulombic cases. The method is also capable of incorporating  multiplicative noise of transport type into the dynamics. It relies in extending functional inequalities of \cite{Serfaty2020, Rosenzweig2020spv, Serfaty2021} to more general interactions, via a new, robust proof that exploits a certain commutator structure. }
\vskip0.3cm



\section{Introduction} \subsection{Main goal and method}
In this article, we consider the first-order mean-field dynamics of interacting particle systems of the form
\begin{equation}
\label{sys}
\begin{cases}
\dot{x}_i =\displaystyle \frac{1}{N}\sum_{{1\leq j\leq N} : j\neq i} \M \nabla\g(x_i-x_j)\\
x_{i}|_{t=0} = x_{i}^0,
\end{cases}
\qquad i\in\{1,\ldots,N\}.
\end{equation}Here, $N\in\N$ is the number of particles, $x_{i}^0 \in \R^d$ are the pairwise distinct initial configurations and $\M$ is a matrix satisfying 
\begin{equation}\label{defM}
\forall \xi \in \R^d, \quad   \langle \M \xi, \xi\rangle \le 0,
\end{equation}
where $\langle \cdot, \cdot \rangle$ is the standard Euclidean scalar product.
Taking $\M  = -\I$ yields  \emph{gradient-flow/ \\ dissipative} dynamics, while taking $\M$ to be an antisymmetric matrix yields {\it Hamiltonian} dynamics. General $\M$ allow for possible mixed flows. We assume that the potential $\g$ is repulsive, so that particles never collide and there is a unique global solution to the system of ODEs \eqref{sys}. The model case we have in mind is $\g$ a Riesz potential depending on a parameter $0\leq s<d$ according to the convention
\begin{equation}
\label{eq:gRiesz}
\g(x) =
\begin{cases}
-\log|x| , & {s=0}\\
\quad |x|^{-s}, & {s>0}.
\end{cases}
\end{equation}
We will in fact consider a slightly more general class of potentials with singularity of the type \eqref{eq:gRiesz}, as explained in the next subsection. Based on the value of the parameter $s$ and its relation to the dimension $d$, we classify the potential $\g$ as \emph{sub-Coulombic} if $0\leq s<d-2$, \emph{Coulombic} if $s=d-2$, and \emph{super-Coulombic} if $d-2<s<d$.

Motivations for studying systems of the form \eqref{sys} are numerous: interacting particles in physics, particle approximation to PDEs, finding equilibrium states of interaction energies, biological and sociological models, large neural networks (see for instance the introduction in \cite{Serfaty2020}). Riesz interactions are particularly interesting for physics and approximation theory, see for instance \cite{DRAW2002, BHS2019}.

By mean-field limit, we mean proving the convergence as $N \to \infty$ of the {\it empirical measure} 
\be\label{munt}
\mu_N^t\coloneqq \frac1N \sum_{i=1}^N \delta_{x_i^t}\ee  associated to a solution $\ux_N^t \coloneqq (x_1^t, \dots, x_N^t)$ of the system \eqref{sys}. If the points $x_i^0$, which themselves depend on $N$, are such that $\mu_N^0$   converges  to some regular measure $\mu^0$, then a formal derivation leads to the expectation that for $t>0$, $\mu_N^t$ converges to the solution of the Cauchy problem with initial data $\mu^0$ for the limiting evolution equation 
\begin{equation}\label{limeq}
\begin{cases}
\partial_t \mu= -\div ((\M \nab \g*\mu) \mu)  & \\
\mu(0)= \mu^0, & \end{cases} \qquad (t,x)\in\R_+\times\R^d
\end{equation}
with $*$ denoting the usual convolution. Proving the convergence of the empirical measure is closely related to proving {\it propagation of molecular chaos} (see \cite{Golse2016ln,HM2014,Jab2014} and references therein) which means showing that if $f_N^0(x_1, \dots, x_N)$ is the initial probability density of seeing particles at positions $(x_1, \dots, x_N)$ and if $f_N^0$ converges to some factorized state  $(\mu^0)^{\otimes N}$, then  the $k$-point marginals $f_{N,k}^t$ converge for all time to $(\mu^t )^{\otimes k}$. Our result implies a convergence of this type as well (see \cref{rem:pc} below).

We note that we restrict to the case $s<d$ because for $s\ge d$, called the {\it hypersingular case}, the potential $\g$ is not integrable near the origin and the convolution $\g* \mu$ no longer makes sense for a nonnegative measure $\mu$. In fact, the expected evolution equation is no longer of the form \eqref{limeq}, see for instance \cite{HSST2020}. This is the reason for restricting ourselves to the {\it potential case} $s<d$. 

Mean-field limits for systems of particles with regular (typically, Lipschitz) interactions have been understood for a long time \cite{Dobrushin1979, Sznitman1991}, with a trajectorial approach consisting in comparing the positions of the particles to the characteristics of the limiting PDE. In contrast, systems with singular interaction are far more challenging and only in recent years have breakthroughs been made to cover the full potential case $s<d$. Various approaches have been put forward to study them. They usually involve finding a good metric to measure the distance of the empirical measure to its evolving expected limit and showing a Gronwall relation on the time evolution of that metric.  The use of the $\infty$-Wasserstein metric has been successful in the sub-Coulombic case $s<d-2$, \cite{Hauray2009,CCH2014}. 
The modulated-energy approach of \cite{Duerinckx2016,Serfaty2020}, itself inspired from \cite{Serfaty2017}, has, in contrast, been successful in the Coulomb and super-Coulombic case $d-2\le s <d$ (in low dimension and the dissipative case in \cite{Duerinckx2016}, and generally in \cite{Serfaty2020}). There, the metric used, called the modulated energy, is a Coulomb/Riesz-based metric, based on the interaction itself, and which can be understood as a renormalized negative-order homogeneous Sobolev norm. To be more precise, the  modulated energy or (squared) ``distance" between $\frac{1}{N}  \sum_{i=1}^N \delta_{x_i} $ and $\mu$  is defined to be
\begin{equation}\label{modenergy} \Fr_N(\ux_N, \mu) \coloneqq \int_{(\R^d)^2 \backslash \triangle} \g(x-y) d\( \frac{1}{N} \sum_{i=1}^N \delta_{x_i} - \mu\) (x) d\( \frac{1}{N} \sum_{i=1}^N \delta_{x_i} - \mu\) (y),
\end{equation}
where we excise the diagonal $\triangle$ in order to remove the infinite self-interaction of each particle. Again, a Gronwall-type inequality can be proved for this modulated energy.

In parallel, work of Jabin and Wang \cite{JW2016, JW2018} put forward a relative-entropy based method which can handle moderately singular interactions, and is particularly tailored for systems with added noise/diffusion. The two approaches (modulated energy and relative entropy) were recently combined in \cite{BJW2019crm, BJW2019edp, BJW2020} into a \emph{modulated free energy} method capable of treating {\it dissipative} singular flows with additive noise. In the course of that combination, these works also  relaxed the assumptions on the interaction (their setting is the torus, which changes the long-range effects though), in particular allowing even for attractive potentials, such as those of Patlak-Keller-Segel-type. 
Note that the one-dimensional case can be handled via a Wasserstein-gradient-flow approach \cite{CFP2012,BO2019} using the fact that in (and only in) dimension 1, the Riesz interaction is convex.

Our present concern is showing that the modulated-energy method functions beyond the Coulombic or super-Coulombic Riesz case. Proving an exact Gronwall relation on the modulated energy requires showing that the derivative of the modulated energy along particle displacement is bounded by the modulated energy itself. This is in essence a functional inequality, which can be proven independently of any dynamics. Plugging into \eqref{modenergy} the solution  $\ux_N^t$ of \eqref{sys} and the solution $\mu^t$ of \eqref{limeq}, direct computation (see \cite[Lemma 2.1]{Serfaty2020}) yields the inequality 
\begin{equation}
\label{eq:MErhs}
\begin{split}
\frac{d}{dt}\Fr_N(\ux_N^t,\mu^t) &\le \int_{(\R^d)^2\setminus\triangle} \paren*{u^t(x)-u^t(y)}\cdot\nabla\g(x-y)d\(\frac{1}{N}\sum_{i=1}^N\d_{x_i^t}-\mu^t\)^{\otimes 2}(x,y),
\end{split}
\end{equation}
where $u^t\coloneqq \M \nab (\g* \mu^t)$. Thus, we need to estimate the right-hand side, which is done by proving the following functional inequality: for any configuration $\ux_N \coloneqq (x_1, \dots, x_N)$, any suitably decaying, bounded probability density $\mu$, and any suitably regular vector field $v : \R^d \to \R^d$, 
\begin{multline} \label{fi}
\int_{(\R^d)^2\backslash \triangle} ( v(x)-v(y)) \cdot \nab \g(x-y) d\( \frac{1}{N} \sum_{i=1}^N \delta_{x_i} - \mu\)^{\otimes 2} (x,y) \le C \Fr_N(\ux_N, \mu)+ C N^{-\alpha}
\end{multline}
for some $C>0$ and $\alpha>0$ depending only on norms of $\mu$ and $v$. The precise statement of the functional inequality \eqref{fi} is given in \cref{prop:rcomm} below.

In the prior works \cite{Duerinckx2016,Serfaty2020} that used the modulated-energy method, the exact Coulomb nature of $\g$ and the fact that the potential $\mathsf{h}\coloneqq \g* \mu$ solves a {\it local equation} in terms of $\mu$ was used crucially as it allowed to rewrite the  left-hand side of \eqref{fi} in terms of the {\it stress-energy tensor} associated to $\mathsf{h}$. Thanks to a dimension extension procedure popularized by Caffarelli and Silvestre \cite{CS2007}, the super-Coulombic case $d-2< s<d$ could also be handled this way by going through a local operator and a stress-energy tensor. This approach was rigid in the exact form of the interaction. Later in \cite{BJW2019edp}, a Fourier-based approach to proving a relation similar to \eqref{fi} was employed and allows, in the case of the torus, to give more flexibility on the interaction, which still needs to have Riesz-type behavior in physical and Fourier space. In fact, the assumptions placed there are similar to the general assumptions imposed in this article, a point on which we elaborate more in \cref{ssec:intropot} below. In \cite{Serfaty2020}, in addition to the stress-tensor structure, a truncation or smearing procedure with a {\it truncation radius depending on the points} and a {\it  monotonicity} with respect to the smearing radii were crucially used.  This allows for a ``renormalization'' of the stress-energy tensor and the energy in order  to account for the removal of the diagonal in \eqref{modenergy} and  \eqref{fi}. 
 
Our main achievement here is to show that the stress-energy tensor structure is not needed to prove \eqref{fi}, nor is a Fourier-based approach. Rather, some underlying {\it commutator structure} to the left-hand side of \eqref{fi} can be exploited for all Riesz interactions (with $s<d$). We write here of commutator structure because one may  view the left-hand side of \eqref{fi}  as a renormalization of the (dualized) commutator
\be \int_{\R^d}\paren*{v\cdot \nab (\g* f) - \g* \( \nab \cdot  (vf)\)}dx,
\ee
or equivalently
\begin{equation}
\ipp*{f, \comm{v}{\frac{\nabla}{(-\Delta)^{\frac{d-s}{2}}}}f} _{L^2},
\end{equation} 
where $$f\coloneqq \frac1N\sum_{i=1}^N \delta_{x_i} -\mu.$$

In fact, in the exact Riesz cases, we could also prove \eqref{fi} via off-the-shelf commutator estimates from the harmonic analysis literature on paraproducts, specifically refinements by Li \cite{Li2019} of the classical fractional Leibnitz rules due to Kato-Ponce \cite{KP1988} and Kenig-Ponce-Vega \cite{KPV1993}. Instead, we provide a fairly direct and essentially physical-space proof based on integration by parts. The Coulomb case  $s=d-2$ is the most delicate. In that case--and only  to obtain the sharp estimate--we appeal to a deep result of Christ and Journ\'e \cite{CJ1987} (see also \cite{SSS2019, Lai2020} for extensions) for commutators of Calder\'{o}n-Zygmund operators. Following the terminology of Christ and Journ\'{e}, such expressions are known as \emph{Calder\'o{n} $d$-commutators}, generalizations of Calder\'{o}n's classic one-dimensional commutator \cite{Calderon1980}. With a little work, the structure in the left-hand side of \eqref{fi} allows for a reduction to estimating such a Calder\'{o}n $d$-commutator. We mention that such commutator expressions often appear in studying fluid problems. See for example \cite{Marchand2008, CCCGW2012, HSSS2018}.

The truncation or smearing procedure of \cite{Serfaty2020} is still needed and is done in a slightly different way depending on whether one is in the sub-Coulombic/Coulombic or the super-Coulombic case. In the sub-Coulombic/Coulombic case, the nature of the interaction makes $\g$ superharmonic, which makes the smearing monotone in the right way. The way we handle super-Coulombic interactions is then by noticing that a subharmonic Riesz potential becomes superharmonic when viewed in a larger dimension. We thus appeal again to a dimension-extension procedure (different from that of \cite{Serfaty2020}, however) to perform the  smearing of the Dirac masses  without other modification of the method. This thus provides a  unified treatment of all Riesz interactions. Since the stress-energy tensor structure is no longer used, we do not need the interaction to be the kernel of a local operator and to be exactly an inverse power. Consequently, we can relax the assumptions in a way similar to \cite{BJW2019edp}. 
 
Our approach is, in principle, robust to the inclusion of multiplicative noise of transport type. Treating the case with such noise was done in \cite{Rosenzweig2020spv} for the two-dimensional Coulomb case, with the proof generalizing to higher dimensions. Due to the nonzero quadratic variation of the noise, it involves proving another functional inequality (see \cref{prop:commso} below), which is exactly a second-order version of \eqref{fi}. In \cite{Rosenzweig2020spv}, the Coulomb nature of the interaction was strongly used. But with the commutator-based approach and the new insights of this article, it is then not too difficult to prove the desired functional inequality in all the Riesz cases we are considering. With this second-order functional inequality and assuming the existence of suitably regular solutions to the limiting stochastic PDEs (SPDEs) to justify the computations, one can follow the road map of \cite{Rosenzweig2020spv} to bound the expected magnitude of the modulated energy. Such a bound, of course, implies convergence in law of the empirical measure to the solution of the limiting SPDE.

\vspace{0.3cm}

\textit{Added in proof:} While our manuscript was under review, a new approach to mean-field limits of first- and second-order systems with additive noise was introduced by Bresch \emph{et al.} \cite{BJS2022}. This approach is based on estimates for the BBGKY hierarchy which take advantage of the diffusion to bound weighted $L^p$ norms of the marginals, and has similarities with the hierarchy approach of Lacker \cite{Lacker2021}. In particular, the work \cite{BJS2022} manages to prove the mean-field limit for the repulsive Vlasov-Fokker-Planck equation with Coulomb potential in dimension 2 and improves upon the prior work of Jabin-Wang \cite{JW2018}.

\subsection{Assumptions and statement of main results}
We now state the precise assumptions for the class of interaction potentials we consider. In the statements below and throughout this article, the notation $\indic_{(\cdot)}$ denotes the indicator function for the condition $(\cdot)$.

For $d\geq 2$ and $0\leq s\leq d-2$, we assume the following:
\be\label{ass0} \g(x) = \g(-x) \quad {\forall x\in\R^{d}\setminus\{0\}} \ee
\be \label{ass1a} \lim_{x\rightarrow 0} \g(x) = +\infty \ee
\be \label{ass1} \text{$\exists r_0>0$ such that}\quad \Delta \g  \leq 0 \quad {\text{in $B(0,r_0)$}}\footnotemark \ee
\footnotetext{Here, we mean $\g$ is superharmonic in $B(0,r_0)$ in the sense of distributions, which implies that $\D\g(x)\leq 0$ for almost every $x\in B(0,r_0)$.}
\be \label{ass2}\forall k \geq 0,\qquad 
|\nabla^{\otimes k} \g(x)|\le C\paren*{\frac{1}{|x|^{s+k}} + |\log |x||\indic_{s=k=0}}  \quad \forall {x\in\R^{d}\setminus\{0\}}\ee
\be\label{ass2b}
|x||\nabla\g(x)| + |x|^2|\nabla^{\otimes 2}\g(x)| \leq C\g(x) \quad \forall x\in B(0,r_0){\setminus\{0\}} \ee
\be\label{ass3}
\frac{C_1} {|\xi|^{ d-s}}  \le  \hat \g(\xi)  \le \frac{C_2}{|\xi|^{d-s}} \quad {\forall \xi\in\R^{d}\setminus\{0\}}\\
\ee
and 
\be \label{ass3a} \text{$\exists c_s<1$ such that} \quad \g(y)<c_s\g(x)\quad \forall x,y\in B(0,r_0)\setminus\{0\} \ \text{with} \ |y|\geq 2|x|.\ee
In the cases $s=d-2k\geq 0$, for some positive integer $k$, we also assume that the $(\R^d)^{\otimes (2k+2)}$-valued kernel
\be \label{ass3'} \mathsf{k}(x-y) \coloneqq  (x-y)\otimes \nabla^{\otimes (2k+1)}\g(x-y) \ee
is associated to a Calder\'{o}n-Zygmund operator in $\R^{d}$.\footnote{Sufficient and necessary conditions for this Calder\'{o}n-Zygmund property are explained in \cite[Section 5.4]{Grafakos2014c}. The reader may check that this condition is satisfied if $\g$ is the genuine Coulomb potential and more generally if $\g$  satisfies $|\nabla^{\otimes k}\hat \g(\xi)|\leq \frac{C}{|\xi|^{d+k-s}}~~\forall \xi\in\R^{d}\setminus\{0\}$.} For $s=0$, we replace the assumption \eqref{ass3a} with
\begin{equation} \label{ass3a'}
\text{$\exists {c_0>0}$ such that} \quad \g(x) - \g(y) \geq {c_0} \quad \forall x,y\in B(0,r_0)\setminus\{0\} \ \text{with} \ |y|\geq 2|x|.
\end{equation}

If $\max\{d-2,0\}<s<d$ or $d=1$ and $s=0$, then we expect the local superharmonicity assumption \eqref{ass1} to fail. The idea, though, is that $\g$ becomes superharmonic when viewing it as an interaction potential in an extended space of dimension $d+m$.  Indeed, $\Delta |X|^{-s}= -s (d+m-s-2) |X|^{-s-2}$ in $\mathbb{R}^{d+m}$, which becomes $\leq 0$ if $m$ is large enough.
We will denote points in $\R^{d+m}$ by $X= (x,z)$ with $x\in \R^d$ and $z\in \R^m$. Our assumption for the case $\max\{d-2,0\}<s<d$ or $d=1$ and $s=0$ is thus that there exists an integer $m\ge 0$, $\G:\R^{d+m}\rightarrow\R$, and constants $C_1, C_2, C>0$ depending only on $d,s,m$, such that
\be\label{ass4} \G(x,0)= \g(x) \quad \forall (x,0)\in \R^{d+m}\ee
\be \label{ass5} \G(X)=\G(-X) \quad {\forall X\in \R^{d+m}\setminus\{0\}} \ee
\be \label{ass6a} \lim_{X\rightarrow 0}\G(X) = +\infty \ee
\be \label{ass6} \text{$\exists r_0>0$ such that} \quad \Delta \G \leq 0 \quad \text{in $B(0, r_0) \subset \R^{d+m}$} \ee
\be \label{ass7}
\forall k\geq 0, \quad|\nabla^{\otimes k} \G(X)|\le C\paren*{\frac1{|X|^{s+k}} + |\log |X||\indic_{s=k=0 \wedge d=1}} \quad \forall {X\in\R^{d+m}\setminus\{0\}} \ee
\be\label{ass7b} 
|X||\nabla\G (X)| + |X|^2 |\nabla^{\otimes 2} \G(X)|\le C\G(X) \quad \forall X\in B(0,r_0)\setminus\{0\}\ee
\be\label{ass8}
 \frac{C_1} {|{\Xi}|^{d+m-s}}  \le  \hat\G({\Xi})  \le \frac{C_2}{|{\Xi}|^{d+m-s}} \quad \forall \Xi\in\R^{d+m}\setminus\{0\}.
\ee
We also add
\begin{equation}\label{ass9}
\begin{cases}
\text{$\exists c_s<1$ such that} \enspace \G(Y) <c_s \G(X)  \quad \forall X,Y\in B(0,r_0) \setminus\{0\} \ \text{with} \ |Y|\geq 2|X|, & {s>0} \\
\text{$\exists {c_0}>0$ such that} \enspace \G(X)-\G(Y)\geq c_0 \quad \forall X,Y\in B(0,r_0)\setminus\{0\} \ \text{with} \ |Y|\geq 2|X|, & {s=0}.
\end{cases}
\end{equation}
If $s=d-1$ and $m=1$, then we also assume that the $(\R^d)^{\otimes 4}$-valued kernel
\begin{equation}\label{ass9'}
\K(X-Y)\coloneqq (X-Y) \otimes \nabla^{\otimes 3}\G(X-Y), \qquad \forall X\neq Y
\end{equation}
is associated to a Calder\'{o}n-Zygmund operator in $\R^{d+m}$. We can simultaneously treat the cases $0\leq s\leq d-2$ and $\max\{d-2,0\}<s<d$ or $d=1$ and $s=0$ by setting $m=0$ if $0\leq s\leq d-2$. If $\g$ satisfies the assumptions \eqref{ass0} -- \eqref{ass3a'} in the case $0\leq s\leq d-2$ or \eqref{ass4} -- \eqref{ass9'} in the case $d-2<s<d$, then we call $\g$ an \emph{admissible potential}.

The main result of this article is the following functional inequality for the modulated energy (cf. \cite[Theorem 1]{Serfaty2020}, \cite[Theorem 3.2]{BJW2019edp}).

\begin{theorem}
\label{thm:main}
Let $\g$ be an admissible potential. Assume the equation \eqref{limeq} admits a solution $\mu \in L^\infty([0,T],\P(\R^d)\cap L^\infty(\R^d))$, for some $T>0$, such that 
\begin{multline}\label{nmut}
\|\Dm^{s-d}\nabla\mu\|_{L^\infty([0,T], L^\infty)}\indic_{s\geq d-1} + \|\nabla^{\otimes 2}\g\ast\mu\|_{L^\infty([0,T], L^\infty)}\indic_{s=d-2} \\
+ \Big(\|\nabla^{\otimes 2}\g\ast\mu\|_{L^\infty([0,T], C^{0,\al})}+ \|\nabla\g\ast\mu\|_{L^\infty([0,T], W^{\frac{d(d+m-s)+2m}{2(d+m)}, \frac{2(d+m)}{d+m-2-s}})}\Big)\indic_{s>d-2}< \infty.
\end{multline}
If $s=0$, then also assume that $\int_{\R^d}\log(1+|x|)d\mu^t(x) < \infty$ for all $t\in [0,T]$. Here, $\Dm=(-\D)^{1/2}$, $m$ is the dimension of the extension coordinate space $\R^m$ for $\g$, and $0<\al\leq 1$.

Let $\ux_N$ solve \eqref{sys}. Then there exist positive constants  $C_3, C_4$, depending only on $|\M|$, $\|\mu^0\|_{L^\infty}$ and the norms of $\mu$ controlled by \eqref{nmut} and on the potential only through the constants in assumptions \eqref{ass1}--\eqref{ass9'}, and an exponent $\beta>0$, depending only on $d,s$,  such that
\be\label{distcoul}
F_N(\ux_N^t , \mu^t)
 \le \( F_N(\ux_N^0, \mu^0)+ C_3 N^{-\beta}\) e^{C_4 t} \qquad \forall t\in [0,T].\ee
In particular, using the notation \eqref{munt},  if $\mu_N^0 \rightharpoonup \mu^0$ in the weak-* topology for measures and
\begin{equation}
\lim_{N\to \infty}F_N(\ux_N^0, \mu^0) =0,\end{equation}
then
\begin{equation}
\mu_N^t \rightharpoonup \mu^t \qquad \forall t\in [0,T].
\end{equation}

\end{theorem}
\begin{remark}
The reader may infer from \cref{sec:Gron} the precise dependence of the constants $C_3, C_4$ above on the norms of the solution $\mu$ as well as the size of the exponent $\beta$. We have omitted this precise dependence from the statement of \cref{thm:main} so as to present the result in the most accessible terms.

We have not attempted to optimize the regularity/integrability assumptions on the solution $\mu$. But we expect that one can proceed in the direction of the works \cite{Rosenzweig2020pvmf, Rosenzweig2020coumf} to only require that $\mu$ belongs to a function space which is critical or borderline with respect to the scaling of the equation. For more discussion on the existence of (necessarily unique) solutions satisfying the condition \eqref{nmut}, we refer the reader to \cite[Subsection 1.3]{Serfaty2020}.
\end{remark}

\begin{remark}
In the case $s=0$, the condition $\int_{\R^d}\log(1+|x|)d\mu(x)<\infty$ is to ensure that the convolution $\g\ast\mu$ is a well-defined function. This assumption is purely qualitative, in the sense that none of our estimates will depend on it. One may check through a Gronwall argument that if the initial datum $\mu^0$ satisfies this condition, then $\mu^t$ also satisfies this condition for all $t\in[0,T]$.
\end{remark}

\begin{remark}
Sufficient conditions for the initial modulated energy $\Fr_N(\ux_N^0,\mu^0)$ to tend vanish as $N\rightarrow\infty$ are explained in \cite[Remark 1.2(c)]{Duerinckx2016}. In particular, convergence of the energy for \eqref{sys} to the energy of \eqref{limeq} together with weak-* convergence of the initial empirical measure to the initial datum $\mu^0$ suffice. 
\end{remark}

\begin{remark}
The weak-* convergence of $\mu_N^t$ to $\mu^t$ is quantitative, as shown in \cref{2.4}. In fact, one can also express the convergence in terms of a suitably negative-order Sobolev norm using \cref{2.4}.
\end{remark}

\begin{remark}\label{rem:pc}
As is well-known, \cref{thm:main} implies propagation of chaos for the marginals of the system \eqref{sys} with initial data $\ux_N^0$ randomly chosen according to the law $(\mu^0)^{\otimes N}$. See \cite[Remark 3.7]{Serfaty2020} for details.
\end{remark}

\subsection{Allowable potentials and further extensions}\label{ssec:intropot}
Assumptions \eqref{ass4} - \eqref{ass9'} allow for a large class of interaction potentials beyond the standard Riesz or log cases. Indeed, the simplest, nontrivial extension is the Lenard-Jones-type potential
\begin{equation}
\g(x) = |x|^{-s} + \varphi_1(x)|x|^{-s_1} + \cdots + \varphi_n(x)|x|^{-s_n},
\end{equation}
where $0<s<d$, $0\leq s_n<\cdots<s_1 < s$, and $\varphi_1,\ldots,\varphi_n$ are smooth, compactly supported functions--not necessarily nonnegative. It is straightforward to check that assumptions \eqref{ass4} - \eqref{ass9'} are satisfied.

In the gradient-flow case, we can even consider genuine Lenard-Jones potentials, where $\varphi_1,\ldots,\varphi_n$ are identically equal to real-valued constants $a_1,\ldots,a_n$. Indeed, let $\chi$ be a smooth, radial bump function adapted to the region $|\xi|\geq 2$. Introducing a small parameter $\kappa>0$, we split
\begin{equation}
\g = \g_{good} + \g_{bad},
\end{equation}
where
\begin{equation}
\begin{split}
\hat{\g}_{bad}(\xi) &\coloneqq c|\xi|^{s-d} + \chi(\kappa \xi)\paren*{\sum_{j=1}^n \tl{c}_j |\xi|^{s_j-d}}, \\
\hat{\g}_{good}(\xi) &\coloneqq \paren*{1- \chi(\kappa \xi)}\paren*{\sum_{j=1}^n \tl{c}_j |\xi|^{s_j-d}}.
\end{split}
\end{equation}
Here, the constants $c,\tl{c}_1,\ldots,\tl{c}_{n}$ are normalization constants from taking the Fourier transform. If $d-2< s<d$, then we set $\G_{bad}(X) \coloneqq \g_{bad}(|X|)$ and $\G_{good}(X) \coloneqq \g_{good}(|X|)$ for $X\in \R^{d+m}$. Since the second term in the definition of $\hat{\g}_{bad}$ is smooth, its inverse Fourier transform and its derivatives are rapidly decaying. It is straightforward to check that
\begin{equation}
\g_{bad}(x) \geq \frac{1}{|x|^{s}}- C\sum_{j=1}^n \frac{1}{|x|^{s_j}}- C, \qquad \forall x\neq 0,
\end{equation}
where $C>0$ is some constant. In particular, there exists $0<r_0\ll 1$, such that for all $|x|\leq r_0$, 
\begin{equation}
\g_{bad}(x) \geq \frac{1}{2|x|^s}.
\end{equation}
Furthermore,
\begin{equation}
-\D\G_{bad}(X) \geq s(d+m-s-2)|X|^{-2-s} - C\sum_{j=1}^n |X|^{-s_{j}-2} - C >0,
\end{equation}
for $|X|\leq r_0$ (possibly choosing $r_0$ smaller) and $d+m-s-2> 0$. Lastly, it follows from the radial symmetry of $\G_{bad}$ that
\begin{equation}
C_1 |\Xi|^{-d-m+s}\leq |\hat{\G}_{bad}(\Xi)| \leq C_2 |\Xi|^{-d-m+s}, \qquad \forall \Xi\neq 0.
\end{equation}
Since $\hat{\g}_{good}$ has compact support and is locally integrable, it belongs to $C^\infty$. Moreover, direct computation shows that $\nabla\g_{good} \in \dot{H}^{\frac{d-s}{2}}$ if $2(s_n+1)>s$. Thus, we can rewrite the system \eqref{sys} as
\begin{equation}
\begin{cases}
\dot{x}_i = \displaystyle -\frac{1}{N}\sum_{\substack{1\leq j\leq N \\ j\neq i}} \paren*{ \nabla\g_{bad}(x_i-x_j) + \mathsf{F}(x_i-x_j)} \\
x_i(0) = x_i^0,
\end{cases}
\end{equation}
where $\mathsf{F}\coloneqq  \nabla\g_{good}$ is now a smooth added force. The potential $\g_{bad}$ satisfies assumptions \eqref{ass4} -- \eqref{ass9'}, and the contribution of $\mathsf{F}$ can be handled using \cite[Lemma 2.3]{Serfaty2020}. One may modify the preceding argument to allow for Lenard-Jones potentials with log interactions as well.

Given that Bresch \emph{et al.} \cite{BJW2019crm, BJW2019edp, BJW2020} also consider a general class of singular interaction potentials, including those of Riesz type, let us compare and contrast our \cref{thm:main} to their main result. Bresch \emph{et al.} consider only the gradient-flow case of the system \eqref{sys} in the periodic setting with random initial data. So rather than working with particle trajectories, they work with the associated Liouville equation. They also allow for additive noise in the dynamics, but we shall ignore this aspect for the purposes of this discussion. Through a Gronwall-type argument, they prove a functional inequality for the modulated free energy. As previously remarked, this quantity combines the relative entropy of earlier work \cite{JW2016, JW2018} with the modulated energy of \cite{Duerinckx2016, Serfaty2020}. By standard arguments, their functional inequality then implies propagation of chaos. Instead of the smearing and extension procedure used in our article, Bresch \emph{et al.} introduce a clever construction of a regularized potential \cite[Lemma 4.1]{BJW2019edp}, which tames the singularity of the Dirac masses. To replace the stress-energy tensor argument from \cite{Serfaty2020} based on the Caffarelli-Silvestre extension, they substitute the regularized potential $\g_{reg}$ into the right-hand side of inequality \eqref{eq:MErhs}, use some elementary Fourier analysis \cite[Lemma 5.2]{BJW2019edp} to bound the resulting expression, and finally estimate the error from this substitution. Although it is never mentioned in their work, their Fourier-analytic argument is, in fact, an elementary commutator estimate, a perspective that we develop in detail in this paper.

As in our work, the potentials they consider are even, singular at the origin, and have nonnegative Fourier transforms. Given that the torus is compact, they only need to assume that the potential is in $L^p(\T^d)$, for some $p>1$. Similar to assumptions \eqref{ass2}, \eqref{ass7} and \eqref{ass2b}, \eqref{ass7b}, Bresch \emph{et al.} need pointwise control on the potential and its derivatives to sufficiently high order. They also need a doubling condition analogous to \eqref{ass3a}, \eqref{ass3a'}, \eqref{ass9}. A key difference, though, is that no pointwise control on $\hat{\g}(\xi)$ is needed, aside from nonnegativity. Bresch \emph{et al.} only need some control on the the first derivative of $\hat{\g}$. Also, superharmonicity plays no role in their work, in sharp contrast to our work

It is worth mentioning that the relative entropy portion of the functional used by Bresch \emph{et al.} is unnecessary and in fact one can use their regularized potential construction with their elementary commutator estimate in order to treat both the conservative and gradient flow cases of the particle system \eqref{sys}. Moreover, one can work directly with the particle dynamics, as opposed to with the Liouville equation. In some sense, one of the new contributions of our paper is the recognition that this is possible, albeit using different methods. In forthcoming work by the second and third authors \cite{RS2021}, we will also show that this is possible when one includes so-called additive noise in the dynamics, distinct from the multiplicative noise considered in this article (see \cref{ssec:noiseOv} for a comparison); and, in fact, one can obtain bounds on the modulated energy which are \emph{global} in time.

Compared to the work of Bresch \emph{et al.}, the advantages of our approach are that it is not limited to the periodic setting; it is pathwise in the sense that no randomization of the initial data is needed, avoiding the use of modulated free energy in favor of the modulated energy; and, perhaps most importantly, the rate of convergence is explicit. In \cite{BJW2019edp}, for instance, an explicit rate is never presented. If one does attempt to optimize the rate in their work, then one sees that the use of the regularized potential construction \cite[Lemma 4.1]{BJW2019edp} is a bottleneck, and the rate one does obtain is  suboptimal compared to our work. Additionally, the regularity requirements for the solution $\mu$ to \eqref{limeq} of \cite{BJW2019edp} are more severe than our work. The cost of these advantages is the need for stronger pointwise control on the potential and its Fourier transform, in particular we have to worry about decay at infinity in physical space and blow-up near the origin in Fourier space.

Finally, as commented in the previous subsection, our method of proof also works when multiplicative noise is added to the right-hand side of the $N$-body problem \eqref{sys}. The expected limiting equation is now a stochastic PDE, so that the modulated energy is now a stochastic process. By combining It\^o's lemma for the modulated energy with pathwise analysis in the form of first- and second-order commutator estimates, we can perform a similar Gronwall argument for the expectation of the magnitude of the modulated energy. In fact, one can bound moments of arbitrarily large degree. We defer a precise statement of the model under consideration and the new functional inequality until \cref{sec:noise}.

\subsection{Organization of article}
\label{ssec:introOrg}
Before introducing the basic notation of the article, let us briefly comment on the organization of the remaining body of the article. In \cref{sec:ME}, we introduce the smearing procedure and establish properties of the modulated energy following the outline of \cite[Section 3]{Serfaty2020}, but for the more general class of potentials satisfying assumptions \eqref{ass4} - \eqref{ass9}. In \cref{sec:comm}, we present our new commutator estimate \cref{propcommu}, and in \cref{sec:rcomm}, we use the smearing procedure from \cref{sec:ME} to renormalize our commuator estimate, obtaining \cref{prop:rcomm}. Then in \cref{sec:Gron}, we use the analysis of \cref{sec:rcomm}, specifically \cref{prop:rcomm}, together with a delicate optimization of the floating small parameters in order to close the Gronwall argument for the modulated energy. This then completes the proof of \cref{thm:main}. Lastly, in \cref{sec:noise}, we close the paper by showing how to extend our modulated-energy approach to variants of the system \eqref{sys}  with multiplicative noise added to the dynamics using second-order commutator estimates.

\subsection{Notation}
\label{ssec:introN}
We close the introduction by specifying the basic notation used in the body of the article without further comment.

Given nonnegative quantities $A$ and $B$, we write $A\lesssim B$ if there exists a constant $C>0$, independent of $A$ and $B$, such that $A\leq CB$. If $A \lesssim B$ and $B\lesssim A$, we write $A\sim B$. To emphasize the dependence of the constant $C$ on some parameter $p$, we sometimes write $A\lesssim_p B$ or $A\sim_p B$.

We denote the natural numbers excluding zero by $\N$ and including zero by $\N_0$. Similarly, we denote the positive real numbers by $\R_+$. Given $N\in\N$ and points $x_{1,N},\ldots,x_{N,N}$ in some set $X$, we will write $\ux_N$ to denote the $N$-tuple $(x_{1,N},\ldots,x_{N,N})$. Given $x\in\R^d$ and $r>0$, we denote the ball and sphere centered at $x$ of radius $r$ by $B(x,r)$ and $\p B(x,r)$, respectively. Given a function $f$, we denote the support of $f$ by $\supp f$.

We make frequent use of tensor notation. $A\otimes B$ denotes the usual tensor product between $A=(A_{i_1\cdots i_m}^{j_1\cdots j_n})$ and $B=(B_{i_1'\cdots i_{m'}'}^{j_1'\cdots j_{n'}'})$ with components
\begin{equation}
(A\otimes B)_{i_1\cdots i_m i_1'\cdots i_{m'}'}^{j_1\cdots j_nj_1'\cdots j_{n'}'} = A_{i_1\cdots i_m}^{j_1\cdots j_n}B_{i_1'\cdots i_{m'}'}^{j_1'\cdots j_{n'}'}.
\end{equation}
We identify $2$-tensors as matrices and use the notation $:$ to denote the Frobenius inner product. Given a function $f$, $\nabla^{\otimes k}f$ denotes the $(\R^d)^{\otimes k}$-valued field whose components are given by the partial derivatives $\p_{i_1\cdots i_k}f$, for any $1\leq i_1,\ldots,i_k\leq d$.

We denote the space of probability measures on $\R^d$ by $\P(\R^d)$. When $\mu$ is in fact absolutely continuous with respect to Lebesgue measure on $\R^d$, we shall abuse notation by writing $\mu$ for both the measure and its density function. We denote the Banach space of complex-valued continuous functions on $\R^d$ by $C(\R^d)$ equipped with the uniform norm $\|\cdot\|_{\infty}$. More generally, we denote the Banach space of $k$-times continuously differentiable functions with bounded derivatives up to order $k$ by $C^k(\R^d)$ equipped with the natural norm, and we define $C^\infty \coloneqq \bigcap_{k=1}^\infty C^k$. Functions $f$ which belong to $C^{k-1}(\R^d)$ with $\nabla^{\otimes k} f$ $\alpha$-H\"older continuous are denoted by $C^{k,\al}(\R^d)$, which we equip with the usual inhomogeneous norm. We use the subscript $c$ to denote the subspace with compact support. Also, we denote the Schwartz space by $\Sc(\R^d)$ and the Bessel potential space by $W^{s,p}(\R^d)$.

\subsection{Acknowledgments}
Q.H.N. is supported by the Shanghai Tech University startup fund, the National Natural Science Foundation of China (12050410257), and NSFC Fund's staff support. M.R. is supported by the Simons Foundation through the Simons Collaboration on Wave Turbulence. S.S. is supported by NSF grant DMS-2000205 and by the Simons Foundation through the Simons Investigator program.

\section{The modulated energy}
\label{sec:ME}
In this section, we perform a functional-analytic study of the modulated energy $\Fr_N$ for an arbitrary pairwise distinct configuration $\ux_N\in (\R^d)^N$ and probability density $\mu\in L^\infty$ which has a logarithmic-type decay at infinity in the case $s=0$. This last assumption is purely qualitative (our estimates will not depend on it): it serves to ensure that the convolution $\g\ast\mu$ is well-defined.

We start by considering the case where \eqref{ass1} holds, which is satisfied by all the sub-Coulombic Riesz cases $0<s<d-2$ as well as the log case $s=0$ if $d\geq 2$.

\subsection{Smearing}
Since \eqref{ass1} holds, $\g$ is superharmonic in $B(0, r_0)$. For any sufficiently integrable function $f$, we have
\begin{equation}\label{mvp}
\frac{d}{dr} \dashint_{\partial B(x, r)} fd\H^{d-1} = \frac{1}{d|B(0,1)|r^{d-1}}\int_{ B(x,r)} \Delta fdy,
\end{equation}
where $\H^{d-1}$ denotes the $(d-1)$-dimensional Hausdorff measure on $\R^d$. This immediately implies that $\g$ satisfies a mean-value inequality:
\be \forall B(x,r) \subset B(0,r_0), \qquad \label{mvi} \g(x) \ge \dashint_{\partial B(x, r)} \g dy.\ee
We note that this inequality holds despite the singularity of $\g$ at the origin. Indeed, $\Delta \g $ is locally integrable if $0\leq s<d-2$ by assumption \eqref{ass2}. If $s=d-2$, then $\D\g$ is locally integrable away from the origin, in which case \eqref{mvp} holds. Using an approximation argument, \eqref{mvi} also follows.

Next, for any  $0<\eta<\min\{\hal, \frac{r_0}{2}\}$ and $x\in \R^d$, we let $\delta_x^{(\eta)}$ be the uniform probability measure on $\partial B(x, \eta)$ and set
\begin{equation}\label{defgeta}
\g_\eta \coloneqq \g\ast\delta_0^{(\eta)}.
\end{equation}
It follows directly from  \eqref{mvi}  that 
\be \label{bgeta} \g_\eta (x) \leq \g(x) \quad \forall x \in B(0, r_0-\eta) \setminus\{0\}\ee
and from  \eqref{mvp} and \eqref{ass2} that  
\be \label{bdiffgeta} | \g(x)-\g_\eta(x)|  \le \frac{C \eta^2}{|x|^{s+2}}  \quad \forall |x|\ge 2\eta, \ee
where the constant $C>0$ depends on $r_0$.

Observe that  in view of \eqref{bgeta} and \eqref{ass2}, the self-interaction of the smeared point mass $\d_{x_0}^{(\eta)}$ satisfies the relation
\be \label{selfinter}
\int_{(\R^d)^2}\g(x-y)d\d_{x_0}^{(\eta)}(x)d\d_{x_0}^{(\eta)}(y) = \int_{\R^d}\g_\eta(x)d\delta_{0}^{(\eta)}(x)  
\le \int_{\R^d}\g(x)d\delta_{0}^{(\eta)}(x)  
 = \g_\eta(0) \le C \eta^{-s} .\ee
If $\g(x)= |x|^{-s}$ then a direct computation with a rescaling is that 
\begin{equation}
\int_{(\R^d)^2}\g(x-y)d\d_{x_0}^{(\eta)}(x)d\d_{x_0}^{(\eta)}(y) = \eta^{-s} \int_{(\R^d)^2}\g(x-y)d\d_{0}^{(1)}(x)d\d_{0}^{(1)}(y) = \eta^{-s}\g_1(1),
\end{equation}
where we abuse notation by letting $1$ denote any choice of unit vector. Similarly, if $\g(x)=-\log|x|$, then
\begin{align}
\int_{(\R^d)^2}\g(x-y)d\d_{x_0}^{(\eta)}(x)d\d_{x_0}^{(\eta)}(y) &= -\log\eta +\int_{(\R^d)^2}\g(x-y)d\d_0^{(1)}(x)d\d_0^{(1)}(y) \nn\\
&= -\log\eta + \g_1(1).
\end{align}

\subsection{Monotonicity, coerciveness and local interaction control}
\label{ssec:MEme}
By definition \eqref{modenergy} of $F_N(\ux_N,\mu)$ and the relation \eqref{selfinter}, we may write 
\begin{equation}\label{redef}
\begin{split}
\Fr_N(\ux_N, \mu)&=\lim_{\alpha_i\to 0}\paren*{ \int_{(\R^d)^2} \g(x-y) d\( \frac{1}{N} \sum_{i=1}^N \delta_{x_i}^{(\alpha_i)}-\mu\)^{\otimes 2}(x,y) - \frac1{N^2}\sum_{i=1}^N \int_{\R^d} \g_{\alpha_i}d \delta_0^{\alpha_i}}.
\end{split}
\end{equation}
This is a way to express $\Fr_N$ as a renormalized version of the undefined expression
\begin{equation}
\int \g(x-y) d\paren*{\frac{1}{N} \sum_{i=1}^N \delta_{x_i}-\mu}^{\otimes 2},
\end{equation}
where the renormalization occurs via the smearing of the Dirac masses. 

One may then examine the change in the expression in the limit when increasing $\alpha_i$ to fixed smearing radii $\eta_i$. As first noted in \cite{PS2017,LSZ2017}, the expression is almost decreasing in the radii, with errors that can be explicitly controlled. 
In addition, thanks to the Fourier transform bounds \eqref{ass8}, once the diagonal has been reinserted in this way 
we have 
\begin{equation}\label{equivtonorm}
\left\|   \frac{1}{N} \sum_{i=1}^N \delta_{x_i}^{(\eta_i)}-\mu    \right\|_{\dot{H}^{\frac{s-d}{2}}(\R^d)}^2
\sim   \int_{(\R^d)^2} \g(x-y) d\( \frac{1}{N} \sum_{i=1}^N \delta_{x_i}^{(\eta_i)}-\mu\)^{\otimes 2}(x,y),
\end{equation}
thus allowing us to control a true Sobolev distance between $\frac{1}{N} \sum_{i=1}^N \delta_{x_i}^{(\eta_i)}$ and $\mu$.

The next proposition expresses this crucial monotonicity property (cf. \cite[Proposition 3.3]{Serfaty2020}) and shows that  the modulated energy is bounded from below and coercive, and controls the small scale interactions. It covers potentials with singularity up to and including Coulomb.

\begin{proposition}
\label{prop:MElb}
Let $d\geq 2$ and $0\leq s\leq d-2$. Suppose that $\ux_N\in ( \R^d)^N$ is pairwise distinct and $\mu \in \P(\R^d)\cap L^\infty(\R^d)$. In the case $s=0$, also suppose that $\int_{\R^d}\log(1+|x|)d\mu(x)<\infty$. There exists a constant $C$ depending only on the potential $\g$ through assumptions \eqref{ass1} - \eqref{ass3a'}, such that for every set of truncation parameters $0< \eta_1,\ldots, \eta_N<r_0/2$, we have 
\begin{multline}\label{respropMElb}
\frac{1}{N^2}\sum_{\substack{1\leq i\neq j\leq N\\ |x_i-x_j| \le \hal r_0} }  
\paren*{\g(x_j-x_i)-\g_{\eta_i}(x_j-x_i)}_+  + C^{-1} \left\|\frac{1}{N}\sum_{i=1}^N\d_{x_i}^{(\eta_i)}-\mu\right\|_{\dot{H}^{\frac{s-d}{2}}}^2  \\
\le \Fr_N(\ux_N,\mu) + \frac{1}{N^2}\sum_{i=1}^N \g_{\eta_i}(0 ) + \frac{C}{N}\sum_{i=1}^N\eta_i^2\\
+\frac{C}{N}\sum_{i=1}^N\|\mu\|_{L^\infty}\paren*{\eta_i^{d-s} + (\eta_i^d|\log\eta_i|)\indic_{s=0} + (\eta_i^2|\log\eta_i|)\indic_{s=d-2}} ,
\end{multline}
where $r_0$ is the constant in \eqref{ass1}.
\end{proposition}
\begin{proof}
Adding and subtracting $\d_{x_i}^{(\eta_i)}$ and regrouping terms, we find that
\begin{align}
\nn \Fr_N(\ux_N,\mu) &= \int_{(\R^d)^2}\g(x-y)d\(\frac{1}{N}\sum_{i=1}^N\d_{x_i}^{(\eta_i)}-\mu\)^{\otimes 2}(x,y)-\frac{1}{N^2} \sum_{i=1}^N \int_{\R^d} \g_{\eta_i} d\delta_0^{(\eta_i)}  \\
&\nn -\frac{2}{N}\sum_{i=1}^N\int_{\R^d} \( \g(y-x_i)-\g_{\eta_i}(y-x_i)\) d\mu(y) \\
&\label{first} + \frac{1}{N^2}\sum_{1\leq i\neq j\leq N} \int_{\R^d} \(\g(y-x_i)-\g_{\eta_i}(y-x_i)\)d(\d_{x_j}+\d_{x_j}^{(\eta)})(y).
\end{align} 

From inequalities \eqref{bgeta} and \eqref{bdiffgeta}, we see that
\begin{multline}
\label{first2}
\frac{1}{N^2}\sum_{1\leq i\neq j\leq N} \int_{\R^d} \(\g(y-x_i)-\g_{\eta_i}(y-x_i)\) d(\d_{x_j}+\d_{x_j}^{(\eta_j)})(y) \\
>\frac{1}{N^2}\sum_{\substack{1\leq i\neq j\leq N \\ |x_i-x_j|< r_0 - \eta }} \(\g(x_j-x_i)-\g_{\eta_i}(x_j-x_i)\)_{+} - \frac{C}{N}\sum_{i=1}^N \eta_i^2.
\end{multline}
Next, using \eqref{bgeta} and  \eqref{bdiffgeta} again, we find that
\begin{align}\label{intgeta}
\int_{B(0, r_0-\eta)} |\g(x)-\g_\eta(x)|dx &\le 2 \int_{B(0,2 \eta)} \g(x)dx + C \eta^2 \int_{B(0,r_0) \backslash B(0, 2\eta)}   \frac{dx}{|x|^{s+2}} \nn\\
&\le C\paren*{\eta^{d-s} + (\eta^d |\log\eta|)\indic_{s=0} + (\eta^2|\log\eta|)\indic_{s=d-2}}.
\end{align}
Since $\mu$ is a probability density, we deduce from H\"older's inequality that
\begin{multline}\label{first3}
-\frac{2}{N}\sum_{i=1}^N\int_{\R^d} \paren*{ \g(y-x_i)-\g_{\eta_i}(y-x_i)}d\mu(y) \geq   - {\frac{C}{N}\sum_{i=1}^N\eta_i^2} \\
- \frac{C}{N}\|\mu\|_{L^\infty} \sum_{i=1}^N \paren*{\eta_i^{d-s} + (\eta_i^d|\log\eta_i|)\indic_{s=0} + (\eta_i^2|\log\eta_i|)\indic_{s=d-2}}.
\end{multline}
Inserting \eqref{first2} and \eqref{first3} into \eqref{first} and using \eqref{equivtonorm}, the result follows. 
\end{proof}

Since $\g_\eta(0)\le C (\eta^{-s} + |\log\eta|\indic_{s=0})$, the right-hand side of \eqref{respropMElb} is easily optimized by choosing $\eta_i= N^{-\frac{1}{2+s}} \leq N^{-1/d}$. Despite not having a sign, $\Fr_N$ is bounded below: 
\begin{equation}\label{esFN}
\Fr_N(\ux_N,\mu) \ge - C(1+\|\mu\|_{L^\infty})\paren*{N^{-\frac{2}{2+s}} + (N^{-1}\log N) \indic_{s=0} + (N^{-\frac{2}{d}} \log N)\indic_{s=d-2}} ,
\end{equation}
which, by triangle inequality, implies that
\begin{multline}
|\Fr_N(\ux_N,\mu)| \leq \Fr_N(\ux_N,\mu) + 2C(1+\|\mu\|_{L^\infty})\Big(N^{-\frac{2}{2+s}} + (N^{-1}\log N) \indic_{s=0} \\
 + (N^{-\frac{2}{d}}\log N) \indic_{s=d-2}\Big).
\end{multline}

\subsection{Renormalization in extended space}
\label{ssec:MEext}
Similar to before, we define $\delta_X^{(\eta)}$ to be the uniform probability measure on $\partial B(X, \eta) \subset \R^{d+m}$ and define the mollified potential
\begin{equation}
\G_\eta\coloneqq \G\ast\delta_0^{(\eta)}.
\end{equation}
With the observation that $\G$ is superharmonic in $B(0,r_0)\subset\R^{d+m}$, we have the analogue of \eqref{bgeta} and \eqref{bdiffgeta}:
\begin{equation} \label{bgetaext}
\G_\eta (X) \leq \G(X) \qquad \forall X \in B(0, r_0-\eta) \setminus\{0\}
\end{equation}
and, from \eqref{mvp} and \eqref{ass2},  
\be \label{bdiffgetaext} | \G(X)-\G_\eta(X)|  \le \frac{C \eta^2}{|X|^{s+2}}  \qquad \forall |X|\ge 2\eta, \ee
where $C>0$ depends on $r_0$. From \eqref{bgetaext}, we obtain the analogue of \eqref{selfinter}:
\be \label{selfint}\int_{(\R^{d+m})^2}\G(x-y)d\d_{0}^{(\eta)}(x)d\d_{0}^{(\eta)}(y) \le \G_\eta(0)  \le C\paren*{\eta^{-s} + |\log\eta|\indic_{d=1 \wedge s=0}} .\ee

We may now repeat all of the analysis of \cref{ssec:MEme} in the extended space $\R^{d+m}$.  The modulated energy is still defined by \eqref{modenergy}. The only change is in the mollification procedure which happens in extended space. Analogous to the lower bound of \cref{prop:MElb}, we have a lower bound for the modulated energy in the super-Coulombic case via this mollification procedure.

\begin{proposition}
\label{prop:extMElb}
Let $d\geq 1$ and $\max\{d-2,0\}< s<d$ or $d=1$ and $s=0$. Suppose that $\ux_N\in ( \R^d)^N$ is a pairwise distinct configuration and $\mu\in \P(\R^d)\cap L^\infty(\R^d)$. If $s=0$, also assume that $\int_{\R^d}\log(1+|x|)d\mu(x)<\infty$. There exists a constant $C$ depending on $s,d$ and on $\g$ only through assumptions \eqref{ass6} -- \eqref{ass9}, such that for every $\eta_i <\min\{\hal, \hal r_0\}$, we have 
\begin{multline}\label{respropMElbext}
\frac{1}{N^2}\sum_{\substack{1\leq i\neq j\leq N\\ |x_i-x_j| \le \hal r_0} }  
\paren*{\g(x_j-x_i)-\G_{\eta_i}(x_j-x_i,0)}_+   +  C^{-1} \left\|\frac{1}{N}\sum_{i=1}^N\d_{X_j}^{(\eta_i)}-\tilde \mu\right\|_{\dot{H}^{\frac{s-d-m}{2}}  (\R^{d+m})}^2  \\
\le \Fr_N(\ux_N,\mu) + \frac{1}{N^2}\sum_{i=1}^N \G_{\eta_i}(0) +\frac{C}{N}\sum_{i=1}^N\paren*{\|\mu\|_{L^\infty}\paren*{\eta_i^{d-s} + (\eta_i^d|\log\eta_i|)\indic_{s=0}} +\eta_i^2},
\end{multline}
where $r_0$ is the constant in \eqref{ass6}, $X_i\coloneqq (x_i,0)$, and
\begin{equation}\label{tildemu}
\tl{\mu}(X) \coloneqq \mu(x)\d_{\R^d\times\{0\}}(X).
\end{equation}
\end{proposition}
\begin{proof}
Following the proof of \cref{prop:MElb}, we write
\begin{equation}\label{firstext}
\begin{split}
\Fr_N(\ux_N,\mu) &=\int_{(\R^{d+m})^2} \G(X-Y)d\Bigg(\frac{1}{N}\sum_{i=1}^N\d_{X_i}^{(\eta_i)}-\tl{\mu}\Bigg)^{\otimes 2}(X,Y) \\
&\ph-\frac{1}{N^2}\sum_{i=1}^N \int_{(\R^{d+m})^2}\G(X-Y)d(\d_{0}^{(\eta_i)})^{\otimes 2}(X,Y) \\
&\ph -\frac{2}{N}\sum_{i=1}^N \int_{\R^{d+m}}\paren*{\G(Y-X_i)-\G_{\eta_i}(Y-X_i)}d\tl{\mu}(Y) \\
&\ph + \frac{1}{N^2}\sum_{1\leq i\neq j\leq N} \int_{\R^d} \paren*{\G(Y-X_i)-\G_{\eta_i}(Y-X_i)}d({\d}_{X_j} + {\d}_{X_j}^{(\eta_j)})(Y).
\end{split}
\end{equation}
The third  line reduces to an   integral over $\R^d$, identical to those in \eqref{first3}.
Analogously to \eqref{intgeta}, in view of \eqref{bgetaext}, \eqref{ass7} and \eqref{bdiffgetaext} we have 
\begin{align}
\int_{B(0, r_0 -\eta) \subset \R^d} |\g(x)-\G_\eta(x,0)|dx  &\leq  2 \int_{B(0, 2\eta)} \g(x) +  C\eta^2 \int_{B(0, r_0)\backslash B(0, 2\eta)}\frac{dx}{|x|^{2+s}} \nn\\
&\leq C \left(\eta^{d-s}+(\eta^d |\log \eta|)\indic_{s=0}\right).
\end{align}
The third line in \eqref{firstext} is thus bounded below by
\begin{equation}
- \frac{C}{N}\sum_{i=1}^N \paren*{\eta_i^2 + \|\mu\|_{L^\infty}  \left(\eta_i^{d-s}+(\eta_i^d |\log \eta_i|)\indic_{s=0}\right)}
\end{equation}
as in \eqref{first3}. For the fourth line, we use   \eqref{bdiffgetaext}  and  $\d_{X_j}^{(\eta_j)}$ is a positive measure to obtain
\begin{align}
&\frac{1}{N^2}\sum_{1\leq i\neq j\leq N} \int_{\R^{d+m}} \paren*{\G(Y-X_i)-\G_{\eta_i}(Y-X_i)}d({\d}_{X_j} + {\d}_{X_j}^{(\eta_j)})(Y) \nn\\
&\ge \frac{1}{N^2}\sum_{\substack{1\leq i\neq j\leq N\\ |x_i-x_j|\le r_0-\eta}} \int_{\R^{d+m}} \paren*{\G(Y-X_i)-\G_{\eta_i}(Y-X_i)}d\d_{X_j}(Y) \nn - \frac{C}{N}\sum_{i=1}^N \eta_i^2 \\
&\ge \frac{1}{N^2}\sum_{\substack{1\leq i\neq j\leq N\\ |x_i-x_j|\le r_0-\eta}} \paren*{\g(x_j-x_i) - \G_{\eta_j}(X_j-X_i)}_+ - \frac{C}{N} \sum_{i=1}^N\eta_i^2 .
\end{align}
Finally, noting that we have by virtue of \eqref{ass8} that
\begin{equation}
\int_{(\R^{d+m})^2} \G(X-Y)d\(\frac{1}{N}\sum_{i=1}^N\d_{X_i}^{(\eta_i)}-\tl{\mu}\)^{\otimes 2}(X,Y)
\geq C_1\left\|\frac{1}{N}\sum_{i=1}^N\d_{X_i}^{(\eta_i)}-\tl{\mu}\right\|_{\dot{H}^{\frac{s-d-m}{2}}(\R^{d+m})}^2,
\end{equation}
we  arrive  at the inequality in the statement of the proposition.
\end{proof}

\subsection{Control of microscale interactions}
As recognized for the first time by the third author in \cite{Serfaty2020}, an interesting choice of parameters $\eta_i$ is to take each $\eta_i$ equal to a fixed fraction of the distance from $x_i$ to its nearest neighbor(s). This ensures that the balls $B(x_i,\eta_i)$ are pairwise disjoint, while also keeping the $\eta_i$ large enough that the right-hand side error in \eqref{respropMElb} and \eqref{respropMElbext} is kept small. However, in this paper, we are able to dispense with the requirement that the balls be disjoint; so, we can choose the $\eta_i$ to all be equal to some parameter $\eta$ or $\ep$, to be optimized at the end. Inequalities \eqref{respropMElb} and \eqref{respropMElbext} allow us to control the small-scale interactions.

Indeed, applying them to  $\eta_i=3\ep$ we arrive at the following corollary. 
\begin{corollary}
\label{cor:counting}
Let $0\le s<d$. For any pairwise distinct $\ux_N \in (\R^d)^N $ and $\mu \in \P(\R^d)\cap L^\infty(\R^d)$ (with $\int_{\R^d}\log(1+|x|)d\mu(x)<\infty$ if $s=0$), we have for all $0<\ep<\min\{\frac{1}{8}, \frac{r_0}{8}\}$,
\begin{multline}
\label{counting}
\frac{1}{N^2}\sum_{\substack{1 \le i \neq j \le N\\
|x_i-x_j|\le \ep} } \paren*{\g(x_i-x_j)\indic_{s>0} + \indic_{s=0}} \le C\bigg(  \Fr_N(  \ux_N, \mu) +C\Big( \frac{\ep^{-s} + |\log\ep|\indic_{s=0}}{N}  +\ep^2\\
+  \|\mu\|_{L^\infty}\big(\ep^{d-s} + (\ep^d|\log\ep|)\indic_{s=0} + (\ep^2|\log\ep|)\indic_{s=d-2}\big)\Big) \bigg),
\end{multline}
where $C$ depends on $s,d$ and on $\g$ only through the assumptions \eqref{ass1} -- \eqref{ass9}. 
\end{corollary}
\begin{proof} 
We use assumption \eqref{ass9}, which ensures that for $|x_i-x_j|\le \ep$, we have 
\begin{equation}
{\g(x_i-x_j) - \g_{3\ep}(x_i-x_j)} \ge (1-c_s)\g(x_i-x_j)\indic_{s>0}+c_0\indic_{s=0}.\qedhere
\end{equation}
\end{proof}

\subsection{Coerciveness}
\label{ssec:MEsob}
From \cref{prop:MElb} or \cref{prop:extMElb}, we easily deduce that the modulated energy $\Fr_N$ controls $\frac{1}{N}\sum_{i=1}^N \delta_{x_i} - \mu$ in the weak-* topology (cf. \cite[Proposition 3.6]{Serfaty2020}). 

\begin{proposition}\label{2.4}
Let $\phi \in \Sc(\R^d)$ be a test function. For any $\eta>0$, any $\ux_N$ and $\mu$ as in \cref{prop:MElb} or \cref{prop:extMElb}, we have
\begin{multline}
\Bigg|\int_{\R^d} \phi \, d\Bigg( \frac{1}{N}\sum_{i=1}^N \delta_{x_i} - \mu\Bigg) \Bigg|\le \|\phi\|_{C^{0, \alpha}} \eta^{\alpha} +  C \|\phi\|_{\dot{H}^{\frac{d-s}{2}} } \bigg( \Fr_N(\ux_N, \mu)+ {C}
\frac{\eta^{-s} + |\log\eta|\indic_{s=0}}{N} \\
+C \|\mu\|_{L^\infty}\big(\eta^{d-s} + (\eta^d|\log\eta|)\indic_{s=0} + (\eta^2|\log\eta|)\indic_{s=d-2}\big) + \eta^2\bigg)^{1/2},
\end{multline}
where $C$ depends on $s,d,m$ and the potential $\g$ through assumptions \eqref{ass1} - \eqref{ass9}.
\end{proposition}

\begin{proof}
Fix $\phi\in\Sc(\R^d)$, and let $\tl\phi\in\Sc(\R^{d+m})$ be an extension of $\phi$ such that
\begin{equation}
\|\tl\phi\|_{\dot{H}^{\frac{d+m-s}{2}}(\R^{d+m})} \lesssim \|\phi\|_{\dot{H}^{\frac{d-s}{2}}(\R^d)},
\qquad
\|\tl{\phi}\|_{C^{0,\alpha}(\R^{d+m})} \lesssim \|\phi\|_{C^{0,\alpha}(\R^d)}.
\end{equation}
For instance, take $\tl\phi(x,z) \coloneqq (e^{|z|\D_{\R^{d}}}\phi)(x)$. We may now write
\begin{multline}
\Bigg|\int_{\R^d} \phi d\Bigg( \frac{1}{N}\sum_{i=1}^N \delta_{x_i} - \mu\Bigg)\Bigg|
\leq \Bigg|\int_{\R^{d+m}} \tilde \phi d\Bigg( \frac{1}{N}\sum_{i=1}^N \delta_{X_i}^{(\eta)} - \tilde \mu\Bigg)\Bigg| +\Bigg|\int_{\R^{d+m}} \tilde \phi   \frac{1}{N}\sum_{i=1}^N d( \delta_{X_i} -\delta_{X_i}^{(\eta)})\Bigg|\\
\leq \|\tilde \phi\|_{\dot{H}^{ \frac{m+d-s}{2}} (\R^{d+m})} \Bigg\| \frac{1}{N}\sum_{i=1}^N \delta_{X_i}^{(\eta)} - \tilde \mu\Bigg\|_{\dot{H}^{\frac{s-d-m}{2}} (\R^{d+m})} + \|\phi\|_{C^{0, \alpha}(\R^d)} \eta^\alpha .
\end{multline}
We conclude after applying the estimate \eqref{respropMElb} or \eqref{respropMElbext}.
\end{proof}

By duality, \cref{2.4} implies that the modulated energy controls inhomogeneous Sobolev norms: for $\gamma>\frac{d+2}{2}$,
\begin{multline}
\Bigg\| \frac{1}{N}\sum_{i=1}^N \delta_{x_i} - \mu\,\,\Bigg\|_{H^{-\gamma}}^2 \le C \Fr_N(\ux_N, \mu) \\
+ C(1+\|\mu\|_{L^\infty})\paren*{N^{\frac{-2}{2+s}} + (N^{-1}\log N)\indic_{s=0} + (N^{-\frac{2}{d}}\log N)\indic_{s=d-2}}.
\end{multline}
Indeed, Sobolev embedding implies that $H^{\gamma}(\R^d)\subset C^1(\R^d)$, from which the preceding bound follows. Further refinements of such norm control, such as in terms of Besov spaces, are possible. For more details, we refer to \cite[Proposition 3.10]{Rosenzweig2020pvmf} and \cite[Proposition 3.8]{Rosenzweig2020coumf}.

\section{Commutator estimates}
\label{sec:comm}
In this section, we present commutator estimates without renormalization, which is the most technical component of the paper. This is a purely space-based functional inequality (cf. \cite[Lemma 5.1]{BJW2019edp}) which dispenses with the stress-energy tensor structure of \cite{Duerinckx2016, Serfaty2020}.

\begin{proposition}
\label{propcommu}
Let $0\leq s<d$. Let $\g \in C^\infty(\R^d\setminus\{0\})$ such that $\g(x)=\g(-x)$ in $B(0,r_0)$, and\footnote{The hypothesis \eqref{eq:commgk} is only necessary for all $k$ up to some integer depending only on $d,s$, which may be determined from examination of the proof of the proposition.}
\begin{equation}
\label{eq:commgk}
\forall k\geq 1, \qquad |\nabla^{\otimes k}\g(x)| \lesssim |x|^{-k-s} \quad \forall x\neq 0.
\end{equation}
If $s\leq d-2$, then also assume that
\begin{equation}
\label{eq:commft}
|\hat{\g}(\xi)| \lesssim |\xi|^{s-d} \qquad \forall \xi\neq 0.
\end{equation}
If $s=d-2k$, for some $k\in\N$, also assume that $\g$ is such that the the kernel
\begin{equation}
\label{eq:commCZker}
\k(x-y) \coloneqq (x-y) \otimes \nabla^{\otimes(1+2k)}\g(x-y) \qquad \forall x\neq y
\end{equation}
is associated to a Calder\'{o}n-Zygmund operator. Let $v$ be a Lipschitz continuous vector field from $\R^d$ to $\R^d$.
There exists a constant $C>0$, depending only $s,d$ and the potential $\g$ through \eqref{eq:commgk} -- \eqref{eq:commCZker}, such that for any $f,g\in\Sc(\R^d)$,\footnote{If $s=0$, then implicit is the requirement that the Fourier transforms of $f,g$ vanish sufficiently rapidly at the origin so that the $\dot{H}^{-d/2}$ norm converges. This may be ensured by assuming a log-weighted $L^1$ condition on $f,g$.} we have
\begin{equation}\label{commu2}
\begin{split}
&\Bigg|\int_{(\R^d)^2} (v(x)-v(y) ) \cdot \nab \g(x-y) f(x) g(y)dxdy\Bigg|\\ 
&\leq  C\paren*{\|\nabla v\|_{L^\infty} + \|\Dm^{\frac{d-s}{2}}v\|_{L^{\frac{2d}{d-2-s}}} \indic_{s<d-2}}  \|f\|_{\dot{H}^{\frac{s-d}{2}}}\|g\|_{\dot{H}^{\frac{s-d}{2}}}.
\end{split}
\end{equation}
Consequently, the integral in the left-hand side of \eqref{commu2} extends to a bounded bilinear form $B_v(\cdot,\cdot)$ on $\dot{H}^{\frac{s-d}{2}}(\R^d)$ satisfying the bound \eqref{commu2}.
\end{proposition}
\begin{proof}
We divide the proof into the super-Coulombic, Coulombic, and sub-Coulombic cases, beginning with the super-Coulombic. We will not explicitly track the dependence of the implicit constants; but the reader should keep in mind that they depend on the parameters $s,d$ as well as the implicit constants in \eqref{eq:commgk}, \eqref{eq:commft} as well as the size and smoothness condition constants for the kernel \eqref{eq:commCZker} and the $L^2$ operator norm of the associated Calder\'{o}n-Zygmund operator.

\noindent
{\bf The super-Coulombic case}.
Define
\begin{equation}
\k_v(x,y) \coloneqq (v(x)-v(y))\cdot\nabla\g(x-y) \qquad \forall x\neq y.
\end{equation}
By approximation, we may assume that $v\in C^\infty$. Furthermore, by density in $\dot{H}^{\frac{s-d}{2}}$, we may assume 
that $f,g\in\Sc(\R^d)$ with Fourier support away from the origin. Additionally, by dominated convergence, we may assume that $\g$ has compact support. Write
\begin{equation}
f =  \div\nabla\Delta^{-1} f \eqqcolon \div f_1
\end{equation}
and similarly for $g$. Note that $f_1=(f_1^i),g_1=(g_1^i)$ are $\R^d$-valued vector fields. Then
\begin{equation}
\int_{(\R^d)^2} \k_v(x,y)f(x)g(y)dxdy = \int_{(\R^d)^2} \k_v(x,y)\div f_1(x)\div g_1(y)dxdy.
\end{equation}
We want to integrate by parts to move the derivatives onto the kernel $\k_v$. But since $\k_v$ is singular along the diagonal $x=y$, we need to exercise some care. To this end, we write the preceding right-hand side as
\begin{equation}
\lim_{\ep\rightarrow 0^+} \int_{|x-y|\geq\ep} \k_v(x,y)\div f_1(x)\div g_1(y)dxdy.
\end{equation}
Integrating by parts in $y$, we find that the preceding integral equals
\begin{equation}
\begin{split}
\int_{|x-y|=\ep} \k_v(x,y)f(x)g_1(y)\cdot\frac{(x-y)}{|x-y|}d\H^{2d-1}(x,y) \\
-\int_{|x-y|\geq\ep}\nabla_y \k_v(x,y)f(x) \cdot g_1(y)dxdy,
\end{split}
\end{equation}
where $\H^{2d-1}$ denotes the $(2d-1)$-dimensional Hausdorff measure on $\R^{2d}$. If $s<d-1$, then we can use the crude kernel estimate based on \eqref{ass7}
\begin{equation}
|\k_v(x,y)| \lesssim \min\(\frac{\|\nabla v\|_{L^\infty}}{|x-y|^{s}}, \frac{\|v\|_{L^\infty}}{|x-y|^{s+1}}\)
\end{equation}
to directly estimate the boundary term by
\begin{equation}
\Bigg|\int_{|x-y|=\ep} \k_v(x,y)f(x)g_1(y)\cdot\frac{(x-y)}{|x-y|}d\H^{2d-1}(x,y)\Bigg| \lesssim \|\nabla v\|_{L^\infty} \|f\|_{L^1} \|g_1\|_{L^\infty} \ep^{d-1-s},
\end{equation}
which tends to zero as $\ep\rightarrow 0$. Otherwise, we subtract and add $g_1(x)$ to write
\begin{equation}
\begin{split}
\int_{|x-y|=\ep}\!\!\!\!\!\!\!\!\!\!\!\!\!\!\k_v(x,y)f(x)g_1(y)\cdot\frac{(x-y)}{|x-y|}d\H^{2d-1}_{(x,y)} 
& = 
\int_{|x-y|=\ep}\!\!\!\!\!\!\!\!\!\!\!\!\!\!\k_v(x,y)f(x)(g_1(y)-g_1(x))\cdot\frac{(x-y)}{|x-y|}d\H^{2d-1}_{(x,y)} \\
&\ph  + 
\int_{|x-y|=\ep}\!\!\!\!\!\!\!\!\!\!\!\!\!\!\k_v(x,y)f(x)g_1(x)\cdot\frac{(x-y)}{|x-y|}d\H^{2d-1}_{(x,y)}.
\end{split}
\end{equation}
In the sequel, it will be convenient to sometimes write in coordinates and employ the convention of Einstein summation. Writing $z=x-y$, so that by Taylor's theorem,
\begin{equation}
\k_v(x,y) = -z^{i_1}\p_{i_1} v(x) \cdot \nabla\g(z) + O(\|\nabla^{\otimes 2}v\|_{L^\infty} |z|^{1-s}),
\end{equation}
we see that
\begin{multline}
\int_{|x-y|=\ep} \k_v(x,y)f(x)g_1(x)\cdot\frac{(x-y)}{|x-y|}d\H^{2d-1}(x,y) =O(\|\nabla^{\otimes 2} v\|_{L^\infty}\|fg_1\|_{L^1}\ep^{d-s})\\
+ \int_{\R^d}\int_{\p B(0,\ep)} (z^{i_1}\p_{i_1}v(x)\cdot \nabla\g(z)) f(x) \paren*{g_1(x)\cdot \frac{z}{|z|}} d\H^{d-1}(z)dx. 
\end{multline}
Note that the first term in the right-hand side vanishes as $\ep\rightarrow 0$ since $s<d$ and the second term vanishes by assumption \eqref{ass5} and the invariance of the measure under $z\mapsto -z$.	

For the non-boundary term, we integrate by parts now in $x$ to obtain
\begin{multline}
-\int_{|x-y|\geq\ep}\nabla_y \k_v(x,y)\div f_1(x) \cdot g_1(y)dxdy =  \int_{|x-y|\geq\ep} \nabla_x\nabla_y \k_v(x,y) : (f_1(x)\otimes g_1(y))dxdy \\
+\int_{|x-y|=\ep} \paren*{\frac{(x-y)}{|x-y|}\otimes \nabla_y \k_v(x,y)} : \paren*{f_1(x)\otimes g_1(y)} d\H^{2d-1}(x,y) .
\end{multline}
We subtract and add $f_1(y)$ to write
\begin{equation}
\begin{split}
&\int_{|x-y|=\ep} \paren*{\frac{(x-y)}{|x-y|}\otimes \nabla_y \k_v(x,y)} : \paren*{f_1(x)\otimes g_1(y)}d\H^{2d-1}(x,y) \\
&=  \int_{|x-y|=\ep} \paren*{\frac{(x-y)}{|x-y|}\otimes \nabla_y \k_v(x,y)} : \paren*{(f_1(x)-f_1(y))\otimes g_1(y)}d\H^{2d-1}(x,y) \\
&\ph +  \int_{|x-y|=\ep} \paren*{\frac{(x-y)}{|x-y|}\otimes \nabla_y \k_v(x,y)} : \paren*{f_1(y)\otimes g_1(y)}d\H^{2d-1}(x,y).
\end{split}
\end{equation}
By Taylor's theorem,
\begin{equation}
\label{eq:commbt}
\begin{split}
&\int_{|x-y|=\ep} \paren*{\frac{(x-y)}{|x-y|}\otimes \nabla_y \k_v(x,y)} : \paren*{(f_1(x)-f_1(y))\otimes g_1(y)}d\H^{2d-1}(x,y) \\
&=\int_{|x-y|=\ep} \paren*{\frac{(x-y)}{|x-y|}\otimes \nabla_y \k_v(x,y)} :\paren*{ (x-y)^{i_1}\p_{i_1}f_1(y)\otimes g_1(y)}d\H^{2d-1}(x,y) \\
&\ph + O(\|v\|_{W^{1,\infty}} \|\nabla^{\otimes 2}f_1\|_{L^\infty}\|g_1\|_{L^1}\ep^{d-s}).
\end{split}
\end{equation}
Proceeding similarly as above, we have that
\begin{equation}
\label{eq:commbt'}
\begin{split}
&\left|\int_{|x-y|=\ep} \paren*{\frac{(x-y)}{|x-y|}\otimes \nabla_y \k_v(x,y)} : \paren*{(x-y)^{i_1}\p_{i_1}f_1(y)\otimes g_1(y)}d\H^{2d-1}(x,y) \right|\\
&\lesssim \| v\|_{W^{2,\infty}} \|\nabla f_1 g_1\|_{L^1} \ep^{d-s}.
\end{split}
\end{equation}
Using that
\begin{equation}
\nabla_x\nabla_y \k_v(x,y) = \nabla_x\nabla_y \k_v(y,x)
\end{equation}
by assumption \eqref{ass5}, we can swap $x$ and $y$ to obtain
\begin{equation}
\begin{split}
&\int_{|x-y|\geq\ep} \nabla_x\nabla_y \k_v(x,y) : \paren*{f_1(x)\otimes g_1(y)}dxdy \\
&= -\frac{1}{2}\int_{|x-y|\geq\ep} \nabla_x\nabla_y \k_v(x,y) : \paren*{(f_1(x)-f_1(y)) \otimes (g_1(x)-g_1(y))}dxdy \\
&\ph + \frac{1}{2}\int_{|x-y|\geq\ep} \nabla_x\nabla_y \k_v(x,y) : \paren*{f_1(x)\otimes g_1(x)}dxdy \\
&\ph + \frac{1}{2}\int_{|x-y|\geq\ep} \nabla_x\nabla_y \k_v(x,y) : \paren*{f_1(y)\otimes g_1(y)}dxdy.
\end{split}
\end{equation}
Note that the last two terms are equal by symmetry and by Stokes theorem,
\begin{equation}
\begin{split}
&\int_{|x-y|\geq\ep} \nabla_x\nabla_y \k_v(x,y) : \paren{f_1(y)\otimes g_1(y)}dxdy \\
&= -\int_{|x-y|=\ep} \paren*{\frac{(x-y)}{|x-y|}\otimes \nabla_y \k_v(x,y)} : \paren*{f_1(y)\otimes g_1(y)}d\H^{2d-1}(x,y).
\end{split}
\end{equation}
After a little bookkeeping, we realize we have shown that
\begin{multline}
\int_{|x-y|\geq\ep} \k_v(x,y) f(x)g(y)dxdy = O(\|\nabla^{\otimes 2} v\|_{L^\infty}\|fg_1\|_{L^1}\ep^{d-s}) \\
+ O(\|v\|_{W^{2,\infty}}\|\nabla f_1 g_1\|_{L^1}\ep^{d-s})  + O(\|v\|_{W^{1,\infty}} \|\nabla^{\otimes 2}f_1\|_{L^\infty}\|g_1\|_{L^1}\ep^{d-s}) \\
 -\frac{1}{2}\int_{|x-y|\geq\ep} \nabla_x\nabla_y \k_v(x,y) : \paren*{(f_1(x)-f_1(y)) \otimes (g_1(x)-g_1(y))}dxdy.
\end{multline}
The first three terms of the right-hand side vanish as $\ep\rightarrow 0$ since $s<d$. For the last term, direct computation reveals
\begin{equation}
\begin{split}
\label{eq:Ksd}
\nabla_x\nabla_y \k_v(x,y) &= -(v(x)-v(y)) \cdot \nabla^{\otimes 3}\g(x-y) - \paren*{\nabla v(y) + \nabla v(x)} : \nabla^{\otimes 2}\g(x-y).
\end{split}
\end{equation}
Using the triangle inequality and mean-value theorem, we find from assumption \eqref{eq:commgk} that
\begin{equation}
|\nabla_x\nabla_y \k_v(x,y)| \lesssim \frac{\|\nabla v\|_{L^\infty}}{|x-y|^{s+2}}.
\end{equation}
This bound together with Cauchy-Schwarz implies that
\begin{equation}
\begin{split}
&\int_{|x-y|\geq\ep} |\nabla_x\nabla_y \k_v(x,y)| |f_1(x)-f_1(y)||g_1(x)-g_1(y)|dxdy \\
&\lesssim \|\nabla v\|_{L^\infty}\paren*{\int_{(\R^d)^2} \frac{|f_1(x)-f_1(y)|^2}{|x-y|^{s+2}}dxdy}^{1/2}\paren*{\int_{(\R^d)^2} \frac{|g_1(x)-g_1(y)|^2}{|x-y|^{s+2}}dxdy}^{1/2}.
\end{split}
\end{equation}
Provided that $d-2<s<d$, the preceding right-hand side is equivalent to (for instance, see \cite[Proposition 3.4]{dNPV2012})
\begin{equation}
\|\nabla v\|_{L^\infty} \|f_1\|_{\dot{H}^{1+\frac{s-d}{2}}} \|g_1\|_{\dot{H}^{1+\frac{s-d}{2}}} \lesssim \|\nabla v\|_{L^\infty} \|f\|_{\dot{H}^{\frac{s-d}{2}}} \|g\|_{\dot{H}^{\frac{s-d}{2}}},
\end{equation}
where the ultimate inequality follows from the definition of $f_1,g_1$. In conclusion,
\begin{equation}
\begin{split}
\lim_{\ep\rightarrow 0} \Bigg|\int_{|x-y|\geq\ep} \k_v(x,y) f(x)g(y)dxdy\Bigg| \lesssim \|\nabla v\|_{L^\infty} \|f\|_{\dot{H}^{\frac{s-d}{2}}} \|g\|_{\dot{H}^{\frac{s-d}{2}}}.
\end{split}
\end{equation}

\noindent
{\bf The Coulombic case}.
We can modify the preceding argument to also work in the Coulombic case $s=d-2$. To see this, we observe from above that
\begin{multline}
\label{eq:Couibp}
\int_{(\R^d)^2}\k_v(x,y)f(x)g(y)dxdy\\
 = \lim_{\ep\rightarrow 0}\Bigg( \int_{|x-y|=\ep} \frac{(x-y)_{i_1}}{|x-y|} \p_{y_{i_2}} \k_v(x,y) f_1^{i_1}(y)g_1^{i_2}(y)d\H^{2d-1}(x,y) \\
+ \int_{|x-y|=\ep} \frac{(x-y)_{i_1}}{|x-y|}\p_{y_{i_2}} \k_v(x,y) (f_1(x)-f_1(y))^{i_1} g_1^{i_2}(y)d\H^{2d-1}(x,y) \\
+ \int_{|x-y|\geq\ep} \p_{x_{i_1}}\p_{y_{i_2}}\k_v(x,y) f_1^{i_1}(x)g_1^{i_2}(y)dxdy\Bigg).
\end{multline}
The second boundary term vanishes as $\ep\rightarrow 0$ by \eqref{eq:commbt} and \eqref{eq:commbt'}. We rewrite the first boundary term as
\begin{equation}
\begin{split}
&\int_{\R^d}  f_1^{i_1}(y) g_1^{i_2}(y)\paren*{\int_{\p B(0,\ep)}\paren*{\frac{z_{i_1}}{|z|} \p_{y_{i_2}} \k_v(y+z,y)} d\H^{d-1}(z)}dy.
\end{split}
\end{equation}
It is straightforward to check that as $\ep\rightarrow 0$, this expression tends to
\begin{equation}
\begin{split}
&-\int_{\R^d} f_1^{i_1}(y) g_1^{i_2}(y) \paren*{\ep^{d}\int_{\p B(0,1)}\frac{z_{i_1}}{|z|} \p_{i_3}v^{i_4}(y) z^{i_3} \p_{i_2i_4}\g(\ep z)) d\H^{d-1}(z)}dy \\
&\ph -\int_{\R^d} f_1^{i_1}(y) g_1^{i_2}(y)  \paren*{\ep^{d-1}\int_{\p B(0,1)}\frac{z_{i_1}}{|z|}  (\p_{i_2} v(y) \cdot\nabla\g(\ep z))d\H^{d-1}(z)}dy.
\end{split}
\end{equation}
Using Cauchy-Schwarz and assumption \eqref{eq:commgk}, we can bound this term directly (up to a constant) by
\begin{equation}
\|\nabla v\|_{L^\infty} \|f_1\|_{L^2}\|g_1\|_{L^2} \lesssim  \|\nabla v\|_{L^\infty}\|f\|_{\dot{H}^{-1}} \|g\|_{\dot{H}^{-1}}.
\end{equation}
For the non-boundary term, we write out $\nabla_x\nabla_y \k_v(x,y)$ to obtain
\begin{equation}
\label{eq:CC}
\begin{split}
&-\int_{|x-y|\geq\ep} \paren*{\p_{i_1}\nabla\g(x-y)\cdot\p_{i_2}v(y) + \p_{i_2}\nabla\g(x-y) \cdot\p_{i_1}v(x)} f_1^{i_1}(x)g_1^{i_2}(y))dxdy  \\
&-\int_{|x-y|\geq\ep} (v(x)-v(y))\cdot\nabla\p_{i_1i_2}\g(x-y) f_1^{i_1}(x) g_1^{i_2}(y))dxdy.
\end{split}
\end{equation}
Observe from assumption \eqref{eq:commft} that the kernel $\nabla^{\otimes 2}\g(x-y)$ defines a $(\R^d)^{\otimes 2}$-valued Fourier multiplier with $L^\infty$ symbol, and therefore by Cauchy-Schwarz,
\begin{multline}
\sup_{\ep>0}\Bigg|\int_{|x-y|\geq\ep} \paren*{\p_{i_1}\nabla\g(x-y)\cdot \p_{i_2}v(y) + \p_{i_2}\nabla\g(x-y)\cdot\p_{i_1}v(x)} f_1^{i_1}(x)g_1^{i_2}(y))dxdy\Bigg| \\
\lesssim \|\nabla v\|_{L^\infty}\|f_1\|_{L^2}\|g_1\|_{L^2} \lesssim \|\nabla v\|_{L^\infty} \|f\|_{\dot{H}^{-1}} \|g\|_{\dot{H}^{-1}}.
\end{multline}
For the remaining term, we use Taylor's theorem to write
\begin{align}
(v(x)-v(y))\cdot\nabla\p_{i_1i_2}\g(x-y) &= \paren*{\int_0^1 \nabla v(y+t(x-y))\cdot (x-y) dt}\cdot \nabla\p_{i_1i_2}\g(x-y) \nn\\
&=\paren*{\int_0^1 \nabla v(y+t(x-y)) dt} \k_{i_1i_2}(x-y),
\end{align}
where the $(\R^d)^{\otimes 2}$-valued kernel $\k_{i_1i_2}$ is defined by
\begin{equation}
\k_{i_1i_2}(x-y)\coloneqq (x-y)\otimes \nabla\p_{i_1i_2}\g(x-y).
\end{equation}
By assumption \eqref{eq:commCZker}, $\k_{i_1i_2}$ is a standard Calder\'o{n}-Zygmund kernel which is associated to a Calder\'{o}n-Zygmund operator (see \cite[Theorem 5.4.1]{Grafakos2014c}). Consequently, the second term in \eqref{eq:CC} is a Calder\'{o}n $d$-commutator, which by the Christ-Journ\'{e} theorem \cite{CJ1987}, is $\lesssim$
\begin{equation}
\|\nabla v\|_{L^\infty} \|f_1\|_{L^2}\|g_1\|_{L^2} \lesssim \|\nabla v\|_{L^\infty} \|f\|_{\dot{H}^{-1}} \|g\|_{\dot{H}^{-1}}.
\end{equation}

\medskip
\noindent
{\bf The sub-Coulombic case}.
Let us now show that the above integration-by-parts argument works for the case $0\leq s<d-2$. First, suppose that $d-4\leq s < d-2$. An examination of the above argument shows that
\begin{equation}
\int_{(\R^d)^2} \k_v(x,y)f(x)g(y)dxdy = \lim_{\ep\rightarrow 0} \int_{|x-y|\geq\ep} \nabla_x\nabla_y \k_v(x,y) : f_1(x)\otimes g_1(y)dxdy.
\end{equation}
Substituting in identity \eqref{eq:Ksd} yields
\begin{multline}\label{eq:sect}
\lim_{\ep\rightarrow 0}\Bigg(-\int_{|x-y|\geq\ep} (v(x)-v(y)) \cdot \nabla\p_{i_1i_2}\g(x-y)f_1^{i_1}(x)g_1^{i_2}(y)dxdy \\ - \int_{|x-y|\geq\ep}\paren*{\p_{i_2}v(y)\cdot\nabla\p_{i_1}\g(x-y) + \p_{i_1} v(x)\cdot\nabla\p_{i_2}\g(x-y)}f_1^{i_1}(x)g_1^{i_2}(y)dxdy\Bigg).
\end{multline}
For the first term, we can repeat the proof above for the super-Coulombic ($d-4<s <d-2$) and Coulombic ($s=d-4$) cases with $\g$ replaced by $\p_{i_1i_2}\g$, and writing $f_1$ as the divergence of a $2$-tensor,
\begin{equation}
f_1^{i_1}= \div\nabla\Delta^{-1}f_1^{i_1} \eqqcolon \div f_2^{i_1},
\end{equation}
and similarly for $g_1$. For the second term in \eqref{eq:sect}, we can pass to the limit $\ep\rightarrow 0$ (note that $\nabla^{\otimes 2}\g$ is locally integrable), so that by Fubini-Tonelli, we now consider
\begin{equation}
\int_{\R^d}(\nabla\p_{i_1}\g\ast f_1^{i_1})(y)\cdot \p_{i_2}v(y)g_1^{i_2}(y)dy + \int_{\R^d}\p_{i_1}v(x)f_1^{i_1}(x)\cdot (\nabla\p_{i_2}\g\ast g_1^{i_2})(x)dx.
\end{equation}
Since $f_1,g_1$ satisfy the same bounds, we can use the $x\leftrightarrow y$ swapping symmetry of the kernel to only consider the first term in the preceding expression. Writing $\mathrm{I} = \Dm^{\frac{d-s}{2}-1}\Dm^{\frac{s-d}{2}+1}$, the Cauchy-Schwarz inequality implies
\begin{multline}
\left|\int_{\R^d}(\nabla\p_{i_1}\g\ast f_1^{i_1})(x) \cdot\p_{i_2} v(x)g_1^{i_2}(x)dx \right| \\
\lesssim \|\Dm^{\frac{d-s}{2}-1}(\nabla\p_{i_1}\g\ast f_1^{i_1})\|_{L^2} \|\Dm^{\frac{s-d}{2}+1}(\p_{i_2} v g_1^{i_2})\|_{L^2}.
\end{multline}
By Plancherel and assumption \eqref{eq:commft}, it then follows that
\begin{equation}\label{eq:fPlan}
\|\Dm^{\frac{d-s}{2}-1}(\nabla\p_{i_1}\g\ast f_1^{i_1})\|_{L^2} \lesssim \|f\|_{\dot{H}^{\frac{s-d}{2}}},
\end{equation}
and by duality,
\begin{align}
\|\Dm^{\frac{s-d}{2}+1}(\p_{i_2}v g_1^{i_2})\|_{L^2} &\leq \sup_{h\in\Sc(\R^d): \|h\|_{L^2}\leq 2} \left|\int_{\R^d}\Dm^{\frac{s-d}{2}+1}h(x) \p_{i_2} v(x)g_1^{i_2}(x)dx\right| \nn\\
&\leq \sup_{h\in\Sc(\R^d): \|h\|_{L^2}\leq 2} \|\Dm^{\frac{d-s}{2}-1}(\Dm^{\frac{s-d}{2}+1}h \p_{i_2}v)\|_{L^2} \|\Dm^{\frac{s-d}{2}+1}g_1^{i_2}\|_{L^2} \nn\\
&\lesssim \|g\|_{\dot{H}^{\frac{s-d}{2}}}\sup_{h\in\Sc(\R^d): \|h\|_{L^2}\leq 2} \|\Dm^{\frac{d-s}{2}-1}(\Dm^{\frac{s-d}{2}+1}h \p_{i_2}v)\|_{L^2}. \label{eq:dualPlan}
\end{align}
By fractional Leibnitz rule (for instance, see \cite[Theorem 7.6.1]{Grafakos2014m})
\begin{align}
\|\Dm^{\frac{d-s}{2}-1}(\Dm^{\frac{s-d}{2}+1}h \p_{i_2}v)\|_{L^2} &\lesssim \|h\|_{L^{p_1}}\|\nabla v\|_{L^{q_1}} + \|\Dm^{\frac{s-d}{2}+1}h\|_{L^{p_2}}\|\Dm^{\frac{d-s}{2}-1}\nabla v\|_{L^{q_2}}
\end{align}
for any $1<p_1,q_1,p_2,q_2\leq\infty$ with $p_i^{-1}+q_i^{-1}=2^{-1}$. We choose $p_1=2$ and $p_2 =\frac{2d}{s+2}$. So by Hardy-Littlewood-Sobolev and the $L^p$ boundedness of Riesz transforms, the right-hand side is $\lesssim$
\begin{equation}\label{eq:dualFL}
\|h\|_{L^2}\paren*{\|\nabla v\|_{L^\infty} + \|\Dm^{\frac{d-s}{2}}v\|_{L^{\frac{2d}{d-2-s}}}}.
\end{equation}
After a little bookkeeping, we arrive at the estimate
\begin{multline}
\left|\int_{\R^d}(\nabla\p_{i_1}\g\ast f_1^{i_1})(x) \cdot \p_{i_2}v(x)g_1^{i_2}(x)dx \right| \\
\lesssim \|f\|_{\dot{H}^{\frac{s-d}{2}}}\|g\|_{\dot{H}^{\frac{s-d}{2}}}\paren*{\|\nabla v\|_{L^\infty} + \|\Dm^{\frac{d-s}{2}}v\|_{L^{\frac{2d}{d-2-s}}}}.
\end{multline}

For the general case $0\leq s<d-2$, we let $k\in\N$ be such that $d-2(k+1) \leq s<d-2k$, and proceed by induction on $k$, using the argument just given above for the induction step. The existence of the form extension $B_v$ follows by the usual density argument.
\end{proof}

\section{Renormalization of the commutator estimate}
\label{sec:rcomm}
In order to apply \cref{propcommu} to $f,g$ which fail to be in $\dot{H}^{\frac{s-d}{2}}$, such as the Dirac mass, we now have to deal with the nondiagonal renormalization which is carried out through the smearing procedure. The main result of this section is the following proposition, which we view as a \emph{renormalized commutator estimate} and which is the workhorse of this article.

\cref{prop:rcomm} stated below should be compared to \cite[Proposition 1.1]{Serfaty2020}, which, as previously discussed, treated the exact Coulomb/super-Coulombic Riesz case, as well as the $\log$ case in dimension 1. That result only required that the vector field $v$ be in $W^{1,\infty}$. In contrast, \cref{prop:rcomm} requires in the super-Coulombic case the stronger assumption that $v\in C^{1,\al}$, for any positive $\al\in (0,1]$, and belong to the Bessel potential space \\ $W^{\frac{d(d+m-s)+2m}{2(d+m)}, \frac{2(d+m)}{d+m-2-s}}$. The reason for this difference is precisely the combination of the commutator estimate with the extension/smearing procedure used in the present article. Namely, we want to apply \cref{propcommu} by using the smearing procedure to renormalize. But this forces us to work in the extended space $\R^{d+m}$ when we apply \cref{propcommu} and then to find a way to control the right-hand side of the proposition in terms of the modulated energy and norms of $v$ and $\mu$, which live in $\R^d$. Consequently, we choose an extension $\tl{v}:\R^{d+m}\rightarrow\R^{d+m}$ of $v$ and apply \cref{propcommu} with $\tl{v}$. It is easy to construct an extension satisfying $\|\tl{v}\|_{W^{1,\infty}} \lesssim \|v\|_{W^{1,\infty}}$ simply by multiplying $v$ by a bump function in the new variable. But such an extension does not satisfy the auxiliary Sobolev condition, which is present (unless $s=d-1$) since we always choose $m$ so that $d+m-s-2\geq 0$. Finding an extension which satisfies both conditions essentially amounts to showing that the trace operator
\begin{equation}
\Tr : W^{1,\infty}(\R^{d+m}) \cap W^{\frac{d+m-s}{2},\frac{2(d+m)}{d+m-2-s}}(\R^{d+m}) \rightarrow W^{1,\infty}(\R^d) \cap W^{\frac{d(d+m-s) + 2m}{2(d+m)},\frac{2(d+m)}{d+m-2-s}}(\R^d)
\end{equation}
is surjective with a bounded right inverse. If one could show such a result, which we do not know how to do, then one could drop the $C^{1,\al}$ assumption.

\begin{proposition}
\label{prop:rcomm}
Let $d\geq 1$, $0\leq s<d$. Let $\ux_N\in (\R^d)^N$ be pairwise distinct, and $\mu \in \P(\R^d)\cap L^\infty(\R^d)$.  If $s=0$, assume that $\int_{\R^d}\log(1+|x|)d\mu(x)<\infty$ and if $s\geq d-1$, assume that \\ $\|\Dm^{s-d}\nabla\mu\|_{L^\infty}<\infty$. Let $v$ be a continuous vector field on $\R^d$ such that
\begin{equation}
\begin{cases}
\|\nabla v\|_{L^\infty} + \|\Dm^{\frac{d-s}{2}}v\|_{L^{\frac{2d}{d-2-s}}}\indic_{s<d-2} < \infty, & {0\leq s\leq d-2} \\
\|v\|_{C^{1,\al}} + \|v\|_{W^{\frac{d(d+m-s)+2m}{2(d+m)}, \frac{2(d+m)}{d+m-2-s}}} < \infty, & {d-2<s<d}
\end{cases}
\end{equation}
for some $0<\al\leq 1$ and where $m$ is the dimension of the extension space. There exists a constant $C$ depending only $d,s,\al,m$ and on the potential $\g$ through assumptions \eqref{ass1} - \eqref{ass9'}, such that
\begin{multline}
\left| \int_{(\R^d)^2\setminus\triangle} (v(x)-v(y))\cdot\nabla\g(x-y)d\paren*{\frac{1}{N}\sum_{i=1}^N\d_{x_i}-\mu}^{\otimes 2}(x,y)\right| \\
\leq C\Bigg(\|\nabla v\|_{L^\infty}+\|\Dm^{\frac{d-s}{2}}v\|_{L^{\frac{2d}{d-2-s}}}\indic_{s<d-2} 
+\paren*{\|v\|_{C^{1,\al}} + \|v\|_{W^{\frac{d(d+m-s)+2m}{2(d+m)}, \frac{2(d+m)}{d+m-2-s}}}}\indic_{s>d-2}\Bigg)\\
\times\Bigg(\Fr_N(\ux_N,\mu) +C(1+\|\mu\|_{L^\infty})\Big(N^{-\frac{\min\{d-s,2\}}{\min\{d,s+2\}(s+1)}} + (N^{-1}\log N) \indic_{s=0} + (N^{-\frac{2}{d(d-1)}}\log N)\indic_{s=d-2}\Big) \\
+ CN^{-\frac{\min\{d+1,s+3\}}{\min\{d,s+2\}(s+1)}}\Big(\|\Dm^{s+1-d}\mu\|_{L^\infty}\indic_{s<d-1} + \|\Dm^{s-d}\nabla\mu\|_{L^\infty}\indic_{s\geq d-1}\Big)\Bigg).
\end{multline}
\end{proposition}
\begin{proof} 
We need to introduce an extension of $\k_v$ to $\R^{d+m}$. By the surjectivity of the trace theorem for Besov and fractional Sobolev spaces (see \cite[Section 2.7.2]{Triebel2010}),\footnote{This is where the $C^{1,\al}$ assumption is used, even though our estimates in the proof will only depend on the Lipschitz seminorm of the extended vector field $\tl{v}$. It is possible to relax $C^{1,\al}$ to the Besov space $B_{\infty,1}^1$, but this is still stronger than $W^{1,\infty}$.} there exists a vector field $\tl{v}$ on $\R^{d+m}$  whose $(m+1)$-th through $(m+d)$-th components vanish identically, $\tl{v}(x,0) = (v(x),0)$, and
\begin{equation}\label{eq:tlvbnd}
\begin{split}
\|\tl{v}\|_{C^{1,\al}(\R^{d+m})} &\lesssim \|v\|_{C^{1,\al}(\R^d)}, \\
\|\tl{v}\|_{W^{\frac{d+m-s}{2}, \frac{2(d+m)}{d+m-2-s}}(\R^{d+m})} &\lesssim \|v\|_{W^{\frac{d(d+m-s)+2m}{2(d+m)}, \frac{2(d+m)}{d+m-2-s}}(\R^d)}.
\end{split}
\end{equation}
With this notation, we define
\begin{equation}
\label{eq:defKv}
\K_v(X,Y) \coloneqq \nabla\G(X-Y) \cdot (\tl{v}(X)-\tl{v}(Y)) \qquad \forall X\neq Y\in \R^{d+m}.
\end{equation}
We note that from  assumption \eqref{ass7} that
\be \label{borneK}
|\K_v(X,Y) |\le {C\|\nab\tl{v}\|_{L^\infty}} \paren*{ |X-Y|^{-s}\indic_{s>0} + \indic_{s=0}}\ee
and 
\be \label{borneDK} 
|\nab_X \K_v(X,Y)|\le  {C\|\nab\tl{v}\|_{L^\infty}} {|X-Y|^{-s-1}}.\ee
For simplicity of notation, we again let $X_i$ denote $(x_i,0)$.

\medskip

Adding and subtracting $\d_{X_i}^{(\eta)}$, for $\eta$ to be chosen shortly, we obtain the decomposition
\begin{equation}
\int_{(\R^d)^2\setminus\triangle} \k_v(x,y)d\Bigg(\frac{1}{N}\sum_{i=1}^N\d_{x_i}-\mu\Bigg)^{\otimes 2}(x,y) = \Te_1+\Te_2+\Te_3,
\end{equation}
where
\begin{align}
\Te_1 &\coloneqq \int_{(\R^{d+m})^2} \K_v(X,Y)d\Bigg(\frac{1}{N}\sum_{i=1}^N\d_{X_i}^{(\eta)}-\tl{\mu}\Bigg)^{\otimes 2}(X,Y),
\end{align}
\begin{align}
\Te_2 &\coloneqq -\frac{2}{N}\sum_{j=1}^N\int_{(\R^{d+m})^2} \K_v(X,Y)d\tl{\mu}(X)d(\d_{X_j}-\d_{X_j}^{(\eta)})(Y),
\end{align}
\begin{align}
\Te_3 &\coloneqq \frac{1}{N^2}\sum_{1\leq i, j\leq N}\int_{(\R^{d+m})^2\setminus\triangle} \K_v(X,Y)d(\d_{X_i}+\d_{X_i}^{(\eta)})(X)d(\d_{X_j}-\d_{X_j}^{(\eta)})(Y).
\end{align}
As before, $\tl{\mu}\coloneqq \mu\d_{\R^d\times\{0\}}$. We proceed to estimate each of the $\Te_j$ individually.

\medskip
\noindent{\bf Estimate for $\Te_1$}.
This term is exactly of the form that can be handled by \cref{propcommu} with
\begin{equation}
f = g=  \frac{1}{N}\sum_{i=1}^N\d_{X_i}^{(\eta)}-\tl{\mu} \in \dot{H}^{\frac{s-d-m}{2}}(\R^{d+m}).
\end{equation}
We deduce from that proposition that 
\be\label{term1}
|\Te_1|\le C\paren*{\|\nabla\tl{v}\|_{L^\infty} + \|\Dm^{\frac{d+m-s}{2}}\tl{v}\|_{L^{\frac{2(d+m)}{d+m-2-s}}} \indic_{s<d+m-2} }  \Bigg\| \frac{1}{N} \sum_{i=1}^N \d_{X_i}^{(\eta)} - \tl\mu \Bigg\|_{\dot{H}^{\frac{s-d-m}{2}}}^2.
\ee
Applying either \cref{prop:MElb} or \cref{prop:extMElb}, we can then bound 
\begin{multline}\label{term1concl}
|\Te_1|\le C \paren*{\|\nabla\tl{v}\|_{L^\infty} + \|\Dm^{\frac{d+m-s}{2}}\tl{v}\|_{L^{\frac{2(d+m)}{d+m-2-s}}} \indic_{s<d+m-2} }  
\\
\times\Bigg( \Fr_N(\ux_N, \mu) + \frac{\G_\eta(0)}{N} + C\paren*{\eta^{d-s} + (\eta^d|\log\eta|)\indic_{s=0} + (\eta^2|\log\eta|)\indic_{s=d-2}}\|\mu\|_{L^\infty} +C\eta^2\Bigg).
\end{multline}

\noindent{\bf Estimate for $\Te_2$}.
Suppose first that $0\leq s<d-1$. We break $\Te_2$ into  a sum over the index $i$ of terms 
\begin{multline}\label{T2sum'}
 \int_{ |X-X_i|< 2\eta}\int_{\R^{d+m}}  \K_v(X,Y) d\( \delta_{X_i} - \delta_{X_i}^{(\eta)}  \) (Y) d\tilde \mu(X) 
\\
+\int_{|X-X_i|\ge 2\eta} \int_{\R^{d+m}}  \K_v(X,Y) d\(  \delta_{X_i} - \delta_{X_i}^{(\eta)}  \) (Y) d\tilde \mu(X).\end{multline}
For the first type of integral, we use that $\d_{X_i}^{(\eta)}$ is a probability measure to rewrite 
\be\label{T2sum}
\int_{\R^{d+m}} \K_v(X,Y) d\( \delta_{X_i} - \delta_{X_i}^{(\eta)}  \) (Y)
= \int_{\R^{d+m}} \(\K_v(X,X_i )- \K_v(X,Y)\)   d \delta_{X_i}^{(\eta)}(Y),
\ee
and unpacking the definition \eqref{eq:defKv} of $\K_v$, we have 
\begin{multline}
\K_v(X,X_i )- \K_v(X,Y)\\
= (\tl{v}(Y) - \tl{v}(X_i)) \cdot \nab \G(X-X_i) + (\tl{v}(X)-\tl{v}(Y)) \cdot ( \nab \G(X-X_i) - \nab \G(X-Y)).
\end{multline}
Using \eqref{ass7}, the mean-value theorem, and triangle inequality
\begin{multline}
\int_{\R^{d+m}} |\K_v(X,X_i )- \K_v(X,Y)|d\delta_{X_i}^{(\eta)}(Y) \\
\leq C\eta\|\nab \tilde v\|_{L^\infty} \bigg(|X-X_i|^{-s-1} +\int_{\R^{d+m}}\paren*{ |X-Y|^{-s-1}}d\delta_{X_i}^{(\eta)}(Y)\bigg)
\end{multline}
for all $|X-X_i| \leq 2\eta$. Integrating the left-hand side with respect $d\tl\mu(X)$ and using Fubini-Tonelli, we find that the first term in \eqref{T2sum'} is controlled by
\begin{equation}
C\eta\|\nab \tilde v\|_{L^\infty}\|\Dm^{s+1-d}\mu\|_{L^\infty},
\end{equation}
where we also use that $|\cdot|^{-s-1}$ is (up to a constant) the convolution kernel of the Fourier multiplier $\Dm^{s+1-d}$ in $\R^d$.  For the second type of contributions we use the mean-value inequality and \eqref{borneDK} to bound 
\begin{align}
&\Bigg| \int_{|X-X_i|\ge 2\eta} \int_{\R^{d+m}}  \K_v(X,Y) d\(  \delta_{X_i} - \delta_{X_i}^{(\eta)}  \) (Y) d\tilde \mu(X)\Bigg| \nn\\
&\leq C \eta\|\nab \tilde v\|_{L^\infty} \int_{|X-X_i|\ge 2\eta}|X-X_i|^{-s-1} d\tilde \mu(X) \nn
\leq C\eta\|\nab \tilde v\|_{L^\infty}{\|\Dm^{s+1-d}\mu\|_{L^\infty}}.
\end{align}

Now if $d-1\leq s<d$, we integrate by parts in $X$ to write
\begin{multline}
\int_{(\R^{d+m})^2} \K_v(X,Y)d\paren*{\d_{X_i}-\d_{X_i}^{(\eta)}}(Y)d\tl\mu(X) \\
= -\int_{(\R^{d+m})^2}\tl\K_v(X,Y)\cdot d\paren*{\d_{X_i}-\d_{X_i}^{(\eta)}}(Y)d\wt{\nabla\mu}(X) \\
-\int_{\R^{d}}\div\tl{v}(x,0)\paren*{\G(x-x_i,0) - \G_\eta(x-x_i,0)}d\mu(x),
\end{multline}
where $\widetilde{\nabla\mu}(X) \coloneqq \nabla\mu(x)\d_{\R^d\times\{0\}}(X)$ and we have defined the vector-valued kernel
\begin{equation}
\tl\K_v(X,Y) \coloneqq \paren*{\tl{v}(X)-\tl{v}(Y)}\G(X-Y).
\end{equation}
Proceeding as in the case $0\leq s<d-1$ with $\mu$ replaced by $\nabla\mu$, the first term in the preceding right-hand side is bounded by
\begin{equation}
C\eta\|\nabla\tl{v}\|_{L^\infty}\|\Dm^{s-d}\nabla\mu\|_{L^\infty}.
\end{equation}
Using that
\begin{equation}
\int_{|x-x_i|\leq 2\eta}\left|\G(x-x_i,0) - \G_\eta(x-x_i,0)\right|dx \leq C\eta^{d-s}
\end{equation}
by assumption \eqref{ass7} and also using \eqref{bdiffgetaext}, we see from H\"older's inequality that the second term is bounded by
$C\|\nabla\tl{v}\|_{L^\infty}\|\mu\|_{L^\infty}\eta^{d-s}$. After a little bookkeeping, we conclude
\begin{multline}
\label{eq:kpT2fin}
|\Te_2| \le C \|\nabla \tl{v}\|_{L^\infty}\Bigg(\eta\|\Dm^{s+1-d}\mu\|_{L^\infty}\indic_{s<d-1}
+ \paren*{\eta\|\Dm^{s-d}\nabla\mu\|_{L^\infty} + \eta^{d-s}\|\mu\|_{L^\infty}}\indic_{s\geq d-1} \Bigg).
\end{multline}

\vspace{.1in}

\noindent{\bf Estimate for $\Te_3$}.
We first remove the self-interaction by observing that
\begin{equation}
\begin{split}
&\sum_{i=1}^N\int_{(\R^{d+m})^2\setminus\triangle}\K_v(X,Y)d(\d_{X_i}+\d_{X_i}^{(\eta)})(X)d(\d_{X_i}-\d_{X_i}^{(\eta)})(Y) \\
&= -\sum_{i=1}^N \int_{(\R^{d+m})^2}\K_v(X,Y)d(\d_{X_i}^{(\eta)})^{\otimes 2}(X,Y),
\end{split}
\end{equation}
so that
\begin{equation}
\begin{split}
\Te_3 &= -\frac{1}{N^2}\sum_{i=1}^N \int_{(\R^{d+m})^2}\K_v(X,Y)d(\d_{X_i}^{(\eta)})^{\otimes 2}(X,Y) \\
&\ph +\frac{1}{N^2}\sum_{1\leq i\neq j\leq N}\int_{(\R^{d+m})^2}\K_v(X,Y)d(\d_{X_i}+\d_{X_i}^{(\eta)})(X)d(\d_{X_j}-\d_{X_j}^{(\eta)})(Y).
\end{split}
\end{equation}
For the first sum, we use assumption \eqref{ass7b} with \eqref{borneK} and \eqref{selfinter} to bound it as 
\begin{align}
\frac{1}{N^2}\Bigg|\sum_{i=1}^N\int_{(\R^{d+m})^2}\K_v(X,Y)d(\d_{X_i}^{(\eta)})^{\otimes 2}(X,Y) \Bigg| &\le \frac{C}{N} \|\nab \tilde v\|_{L^\infty}\paren*{\G_\eta(0)\indic_{s>0} + \indic_{s=0}} \nn\\
&\leq\frac{C}{N} \|\nab \tilde v\|_{L^\infty}\paren*{\eta^{-s}\indic_{s>0} + \indic_{s=0}}. \label{bornet3}
\end{align}

Next, we split the second sum over $1\leq i\neq j\leq N$ into a sum over ``close pairs" for which $|x_i-x_j|< \ep$ and a sum over ``far pairs" for which $|x_i-x_j|\ge \ep$, where $\ep>2\eta$ is a parameter, to be chosen later.
For the close sum, we again use \eqref{ass7b} with the bound \eqref{borneK} to obtain
\begin{align}
&\frac{1}{N^2}\sum_{{1\leq i\neq j\leq N}\atop{|x_i-x_j|< \ep}} \left|\int_{(\R^{d+m})^2}\K_v(X,Y)d(\d_{X_i}+\d_{X_i}^{(\eta)})(X)d(\d_{X_j}-\d_{X_j}^{(\eta)})(Y)\right| \nn\\
& \le   \frac{C\|\nabla \tl{v}\|_{L^\infty}}{N^2} \sum_{{1\leq i\neq j\leq N}\atop{|x_i-x_j|< \ep}}\Bigg(\(
\g(x_i-x_j) + \int \G_\eta(X_i-Y) d\delta_{X_j}^{(\eta)}(Y)\)\indic_{s>0} +\indic_{s=0}\Bigg). \label{eq:LHSleq}
\end{align}
Using \eqref{bgetaext}, we then obtain (provided $\ep < \frac{r_0}{2}$) that the left-hand side of \eqref{eq:LHSleq} is $\leq$
\begin{align}
&C \frac{\|\nabla \tl{v}\|_{L^\infty}}{N^2} \sum_{{1\leq i\neq j\leq N}\atop{|x_i-x_j|< \ep}} \paren*{\g(x_i-x_j)\indic_{s>0} + \indic_{s=0}} \nn\\
&\le C \|\nabla\tl{v}\|_{L^\infty}\Bigg(\Fr_N(\ux_N,\mu) +C\Big(\frac{\ep^{-s}+|\log\ep|\indic_{s=0}}{N} +\ep^2 \nn\\
&\ph +  \|\mu\|_{L^\infty}(\ep^{d-s}+(\ep^d|\log\ep|)\indic_{s=0} + (\ep^2|\log\ep|)\indic_{s=d-2})\Big)\Bigg)\label{eq:far},
\end{align}
where the ultimate inequality follows from application of \cref{cor:counting}. For the sum over far pairs, we first note that by symmetry under swapping $i\leftrightarrow j$,
\begin{equation}
\begin{split}
&\sum_{{1\leq i,j\leq N}\atop {|x_i-x_j|\ge \ep}} \int_{(\R^{d+m})^2}\K_v(X,Y)d(\d_{X_i}+\d_{X_i}^{(\eta)})(X)d(\d_{X_j}-\d_{X_j}^{(\eta)})(Y)\\
&=\sum_{{1\leq i,j\leq N}\atop {|x_i-x_j|\ge \ep}} \int_{(\R^{d+m})^2}\paren*{\K_v(X_i,X_j)-\K_v(X,Y)}d\d_{X_i}^{(\eta)}(X)d\d_{X_j}^{(\eta)}(Y).
\end{split}
\end{equation}
Since $\ep>2 \eta$, using the mean-value inequality and \eqref{borneDK}, we find
\begin{equation}\label{borne2t3}
\frac{1}{N^2}\sum_{{1\leq i\neq j\leq N}\atop{|x_i-x_j|\ge \ep}} \left|\int_{(\R^{d+m})^2}\K_v(X,Y)d(\d_{X_i}+\d_{X_i}^{(\eta)})(X)d(\d_{X_j}-\d_{X_j}^{(\eta)})(Y)\right|   \le \frac{C \eta \|\nab \tilde v\|_{L^\infty}}{\ep^{s+1}}.
\end{equation}
Putting together the estimates \eqref{bornet3}, \eqref{eq:far}, and \eqref{borne2t3}, we obtain
\begin{multline}
\label{eq:kpT3fin}
|\Te_3| \le C  \|\nabla \tl{v}\|_{L^\infty}\Bigg(\Fr_N(\ux_N,\mu) + C\Big(\frac{\eta^{-s} + |\log\eta|\indic_{s=0}}{N}  +\frac{\ep^{-s} + |\log\ep|\indic_{s=0}}{N} + \ep^2 \\
+ \|\mu\|_{L^\infty}(\ep^{d-s}+(\ep^d|\log\ep|)\indic_{s=0} + (\ep^2|\log\ep|)\indic_{s=d-2}) + \eta\ep^{-s-1} \Big) \Bigg).
\end{multline}
\noindent
{\bf Conclusion.}
Combining \eqref{term1concl}, \eqref{eq:kpT2fin}, and \eqref{eq:kpT3fin}, we have shown that
\begin{multline}
\Bigg|\int_{(\R^d)^2\setminus\triangle} \k_v(x,y)d\Bigg(\frac{1}{N}\sum_{i=1}^N\d_{x_i}-\mu\Bigg)^{\otimes 2}(x,y)\Bigg|
\leq C \paren*{\|\nabla \tl{v}\|_{L^\infty} + \|\Dm^{\frac{d+m-s}{2}}\tl{v}\|_{L^{\frac{2(d+m)}{d+m-2-s}}} \indic_{s<d+m-2}}  
\\
\times\Bigg( \Fr_N(\ux_N, \mu) + C\Big(\frac{\eta^{-s}+|\log\eta|\indic_{s=0}}{N} +\eta^2 
+ \|\mu\|_{L^\infty}(\eta^{d-s} +(\eta^d|\log\eta|)\indic_{s=0} + (\eta^2|\log\eta|)\indic_{s=d-2})\Big)\Bigg)
\\
+C  \|\nabla\tl{v}\|_{L^\infty}\Bigg((\eta \|\Dm^{s+1-d}\mu\|_{L^\infty})\indic_{s<d-1} + (\eta\|\Dm^{s-d}\nabla\mu\|_{L^\infty} + \eta^{d-s}\|\mu\|_{L^\infty})\indic_{s\geq d-1} \\
+\frac{\eta^{-s} + |\log\eta|\indic_{s=0}}{N}  +\frac{\ep^{-s} + |\log\ep|\indic_{s=0}}{N} + \ep^2 \\
+ \|\mu\|_{L^\infty}(\ep^{d-s}+(\ep^d|\log\ep|)\indic_{s=0} + (\ep^2|\log\ep|)\indic_{s=d-2}) + \eta\ep^{-s-1} \Bigg).
\end{multline}
We may now optimize over $\eta$ and $\ep$ by taking 
\begin{equation}
\label{eq:rcommopt}
\begin{cases} 
\eta= \ep^{s+2}\ \text{and} \ \ep = N^{-1} , & {s=0}\\
\eta = \ep^{s+3}\ \text{and} \ \ep = N^{-\frac{1}{(s+2)(s+1)}}, & {0<s\leq d-2}\\
\eta= \ep^{d+1} \ \text{and} \ \ep = N^{-\frac{1}{d(s+1)}}, &{d-2<s<d}.
\end{cases}
\end{equation}
Using \eqref{eq:tlvbnd} to control the norms of $\tl{v}$ in terms of norms of $v$, we then find 
\begin{multline}
\label{eq:rcommrhs2}
\Bigg|\int_{(\R^d)^2\setminus\triangle} \k_v(x,y)d\Bigg(\frac{1}{N}\sum_{i=1}^N\d_{x_i}-\mu\Bigg)^{\otimes 2}(x,y)\Bigg| \leq C\Bigg(\|\nabla v\|_{L^\infty}+\|\Dm^{\frac{d-s}{2}}v\|_{L^{\frac{2d}{d-2-s}}}\indic_{s<d-2} + \\
+\paren*{\|v\|_{C^{1,\al}} + \|v\|_{W^{\frac{d(d+m-s)+2m}{2(d+m)}, \frac{2(d+m)}{d+m-2-s}}}}\indic_{s>d-2}\Bigg)\Bigg(\Fr_N(\ux_N,\mu) \\
+C(1+\|\mu\|_{L^\infty})\Big(N^{-\frac{\min\{d-s,2\}}{\min\{d,s+2\}(s+1)}} + (N^{-1}\log N)\indic_{s=0} + (N^{-\frac{2}{d(d-1)}}\log N)\indic_{s=d-2}\Big) \\
+ CN^{-\frac{\min\{d+1,s+3\}}{\min\{d,s+2\}(s+1)}}\Big(\|\Dm^{s+1-d}\mu\|_{L^\infty}\indic_{s<d-1} + \|\Dm^{s-d}\nabla\mu\|_{L^\infty}\indic_{s\geq d-1}\Big)\Bigg).\qedhere
\end{multline}
\end{proof}

\section{Gronwall argument}
\label{sec:Gron}
We now have all the ingredients for our Gronwall argument on the modulated energy. Set $u\coloneqq \M\nabla\g\ast\mu$. Applying the bound of \cref{prop:rcomm} pointwise in time with $v=u^t$ to the right-hand side of inequality \eqref{eq:MErhs} and then integrating with respect to time, we find that for $0\leq t\leq T$,
\begin{multline}
|\Fr_N(\ux_N^t,\mu^t)| \leq |\Fr_N(\ux_{N}^0,\mu^0)| 
\\
+ C\int_0^t\Bigg(\|\nabla u^\tau\|_{L^\infty}+\|\Dm^{\frac{d-s}{2}}u^\tau\|_{L^{\frac{2d}{d-2-s}}}\indic_{s<d-2}
+\paren*{\|u^\tau\|_{C^{1,\al}} + \|u^\tau\|_{W^{\frac{d(d+m-s)+2m}{2(d+m)}, \frac{2(d+m)}{d+m-2-s}}}}\indic_{s>d-2} \Bigg) 
\\
\times\Bigg(\Fr_N(\ux_N^\tau,\mu^\tau)
+C(1+\|\mu^\tau\|_{L^\infty})\Big(N^{-\frac{\min\{d-s,2\}}{\min\{d,s+2\}(s+1)}} + (N^{-1}\log N)\indic_{s=0} + (N^{-\frac{2}{d(d-1)}}\log N)\indic_{s=d-2}\Big) \\
+ CN^{-\frac{\min\{d+1,s+3\}}{\min\{d,s+2\}(s+1)}}\Big(\|\Dm^{s+1-d}\mu^\tau\|_{L^\infty}\indic_{s<d-1} + \|\Dm^{s-d}\nabla\mu^\tau\|_{L^\infty}\indic_{s\geq d-1}\Big)\Bigg)d\tau.
\end{multline}
If $0\leq s<d-2$, then by Sobolev embedding, assumption \eqref{ass3}, and H\"older's inequality,
\begin{equation}
\|\Dm^{\frac{d-s}{2}}u^\tau\|_{L^{\frac{2d}{d-2-s}}} +\|\nab u^\tau\|_{L^\infty} \lesssim 1+\|\mu^\tau\|_{L^\infty}.
\end{equation}
If $0\leq s<d-1$, then it also follows from H\"older's inequality that
\begin{equation}
\|\Dm^{s+1-d}\mu^\tau\|_{L^\infty} \lesssim 1+\|\mu^\tau\|_{L^\infty}.
\end{equation}
Applying the Gronwall-Bellman lemma, we conclude that for all $0\leq t\leq T$,
\begin{multline}
|\Fr_N(\ux_N^t,\mu^t)| \leq A_N(t)\exp\Bigg(C\int_0^t \Big(\|\nabla u^\tau\|_{L^\infty} + \|\Dm^{\frac{d-s}{2}}u^\tau\|_{L^{\frac{2d}{d-2-s}}}\indic_{s<d-2} \\
+ \Big(\|u^\tau\|_{C^{1,\al}}+\|u^\tau\|_{W^{\frac{d(d+m-s)+2m}{2(d+m)}, \frac{2(d+m)}{d+m-2-s}}}\Big)\indic_{s>d-2}\Big)d\tau\Bigg),
\end{multline}
where the time-dependent prefactor $A_N$ is defined by
\begin{multline}\nonumber
A_N(t) \coloneqq  C\int_0^t \Bigg((1+\|\mu^\tau\|_{L^\infty})\Big(N^{-\frac{\min\{d-s,2\}}{\min\{d,s+2\}(s+1)}} + (N^{-1}\log N) \indic_{s=0} + (N^{-\frac{2}{d(d-1)}}\log N)\indic_{s=d-2}\Big) \\
+ CN^{-\frac{\min\{d+1,s+3\}}{\min\{d,s+2\}(s+1)}}\Big(\|\Dm^{s+1-d}\mu^\tau\|_{L^\infty}\indic_{s<d-1} + \|\Dm^{s-d}\nabla\mu^\tau\|_{L^\infty}\indic_{s\geq d-1}\Big)\Bigg)d\tau + |\Fr_N(\ux_{N}^0,\mu^0)|.
\end{multline}
This proves our bound for the evolution of the modulated energy, and hence \cref{thm:main}.

\section{Multiplicative Noise}
\label{sec:noise}
\subsection{Overview}
\label{ssec:noiseOv}
In this last section, we show how our method of proof can in principle be extended to treat the mean-field limit of first-order systems, such as \eqref{sys}, with multiplicative noise added to the dynamics:
\begin{equation}
\label{sysnoise}
\begin{cases}
\dot{x}_i =\displaystyle \frac{1}{N}\sum_{{1\leq j\leq N} : j\neq i} \M \nabla\g(x_i-x_j) + \sum_{k=1}^\infty \sigma_k(x_i)\circ \dot{W}_k\\
x_{i}|_{t=0} = x_{i}^0,
\end{cases}
\qquad i\in\{1,\ldots,N\}.
\end{equation} 
Here, the $\sigma_k$ are $C^\infty$ vector fields on $\R^d$, the $W_k$ are independent real Brownian motions, and $\circ$ denotes that product should be interpreted in the Stratonovich sense.\footnote{When the $\circ$ is not present, the product should be interpreted in the It\^o sense throughout this paper.} The correspond mean-field equation is no longer deterministic, but instead gains a stochastic term:
\begin{equation}
\label{eq:mfnoise}
\begin{cases}
\partial_t \mu= \displaystyle-\div((\M \nab \g*\mu) \mu) - \sum_{k=1}^\infty \div(\sigma_k\mu) \circ\dot{W}_k  & \\
\mu(0)= \mu^0 & 
\end{cases}
\qquad (t,x)\in\R_+\times\R^d.
\end{equation}
The model example of a system of the form \eqref{sysnoise} is the {\it stochastic point vortex model} of Flandoli, Gubinelli, and Priola \cite{FGP2011}, for which $\g$ is the 2D Coulomb potential. The limiting SPDE is now the 2D incompressible Euler vorticity equation with an additional stochastic transport term, the well-posedness of which has been studied by Brze\'{z}niak, Flandoli, and Maurelli \cite{BFM2016}. We emphasize that the noise in \eqref{sysnoise} is very different than the so-called additive noise model,
\begin{equation}
\label{eq:anoise}
\dot{x}_i = \frac{1}{N}\sum_{{1\leq j\leq N} : j\neq i} \M \nabla\g(x_i-x_j) + \sqrt{2\sigma}\dot{\ul{W}}_i,
\end{equation}
where now $\sigma>0$ is a constant and the $\ul{W}_i$ are independent Brownian motions in $\R^d$. Indeed, in equation \eqref{eq:anoise}, the noise is independent for each particle, whereas in equation \eqref{sysnoise}, the spatially dependent noise acting on each particle is the same. This difference clearly manifests itself at the limiting level. As has been rigorously shown by numerous authors \cite{Osada1987pc, Osada1987lp, FHM2014, LY2016, JW2018, LLY2019, BJW2019crm, BJW2019edp, BJW2020}, the empirical measure of the system \eqref{eq:anoise} converges in law to the solution of the deterministic diffusive/viscous equation   
\begin{equation}
\p_t\mu = -\div((\M\nab \g\ast\mu)\mu) + \sigma \D\mu.
\end{equation}

As mentioned in the introduction, the second author \cite{Rosenzweig2020spv} was the first to study the convergence in law of the empirical measure for \eqref{sysnoise} to the solution (if one exists) of the SPDE \eqref{eq:mfnoise}.\footnote{An earlier work of Coghi and Maurelli \cite{CM2020} considered this problem but with an $N$-dependent asymptotically vanishing truncation of the potential in \eqref{sysnoise} to distances much larger than the typical interparticle distance $N^{-1/d}$.} By developing a stochastic extension of the modulated-energy approach of \cite{Serfaty2020} and introducing the commutator perspective adopted in this paper, he proved a quantitative estimate for the convergence in the 2D Coulomb case corresponding to the aforementioned stochastic point vortex model. The proof can be extended to dimensions $d\geq 3$, but, importantly, it very strongly depends on the Coulomb nature of the interaction. At the time of the present article, this is the only such result for singular interactions of which we are aware (see \cite{CF2016} for the case of regular interactions).

The idea of \cite{Rosenzweig2020spv} is to consider the modulated energy $\Fr_N(\ux_N^t,\mu^t)$ as before, noting that it is now a real-valued stochastic process which is almost surely finite. A formal application of It\^o's lemma (in Stratonovich form) to $\Fr_N(\ux_N^t,\mu^t)$ yields the stochastic differential inequality
\begin{multline}
\label{eq:noiseME}
\frac{d}{dt}\Fr_N(\ux_N^t, \mu^t) \leq  \int_{(\R^d)^2\setminus\triangle} \paren*{u^t(x)-u^t(y)}\cdot\nabla\g(x-y)d\Bigg(\frac{1}{N}\sum_{i=1}^N\d_{x_i^t}-\mu^t\Bigg)^{\otimes 2}(x,y) \\
+ \sum_{k=1}^\infty\int_{(\R^d)^2\setminus\triangle} \nabla\g(x-y)\cdot \paren*{\sigma_k(x)-\sigma_k(y)} d\Bigg(\frac{1}{N}\sum_{i=1}^N\d_{x_i^t}-\mu^t\Bigg)^{\otimes 2}(x,y) \circ \dot{W}_k,
\end{multline}
where again $u^t \coloneqq \M\nabla\g\ast\mu^t$. Converting from Stratonovich to It\^o, the second line becomes
\begin{multline}\label{eq:StI}
\sum_{k=1}^\infty \int_{(\R^d)^2\setminus\triangle}\nabla\g(x-y)\cdot\paren*{\sigma_k(x)-\sigma_k(y)} d\Bigg(\frac{1}{N}\sum_{i=1}^N\d_{x_i^t}-\mu^t\Bigg)^{\otimes 2}(x,y)\dot{W}_k\\
+\frac{1}{2}\sum_{k=1}^\infty\int_{(\R^d)^2\setminus\triangle} \nabla\g(x-y) \cdot \paren{\sigma_k\cdot\nabla\sigma_k(x)-\sigma_k\cdot\nabla\sigma_k(y)}d\Bigg(\frac{1}{N}\sum_{i=1}^N\d_{x_i^t}-\mu^t\Bigg)^{\otimes 2}(x,y)\\
+\frac{1}{2}\sum_{k=1}^\infty\int_{(\R^d)^2\setminus\triangle} \nabla^{\otimes 2}\g(x-y) : \paren*{\sigma_k(x)-\sigma_k(y)}^{\otimes 2} d\Bigg(\frac{1}{N}\sum_{i=1}^N\d_{x_i^t}-\mu^t\Bigg)^{\otimes 2}(x,y).
\end{multline}

The first term in \eqref{eq:StI}, which is now an It\^o integral, should have zero expectation, therefore we can ignore it. The second term has the same commutator structure as in the functional inequality \eqref{fi} with $v=\sum_{k=1}^\infty \sigma_k\cdot\nabla\sigma_k$. Consequently, we can treat this expression as a renormalized commutator using \cref{prop:rcomm}. The third term is more complicated. It appears similar to the first two terms, except now there are two derivatives on $\g$ and two symmetrized expressions $\sigma_k(x)-\sigma_k(y)$ which, in principle, should cancel out the singularities induced by the derivatives. However, \cref{prop:rcomm} does not cover such expressions. Accordingly, the goal of this section is to prove its second-order analogue, which is \cref{prop:rcommso} stated below. Again, the idea will be to prove a second-order version of the commutator estimate of \cref{propcommu} using integration by parts and then to obtain a renormalization of this estimate using the smearing procedure.

The preceding computation to arrive at \eqref{eq:noiseME} and \eqref{eq:StI} is completely formal. We cannot directly invoke It\^o's lemma because the potential $\g$ is singular. Also, we have not specified whether there is a solution to the system \eqref{sysnoise} or in what sense the equation \eqref{eq:mfnoise} holds. The $N$-body problem is relatively straightforward to make sense of. Using a truncation of the potential near the origin and a stopping time argument in the spirit of \cite[Section 3]{FGP2011}, one can show there is a unique strong solution. In particular, with probability one, the particles never collide. The SPDE \eqref{eq:mfnoise} is more difficult to make sense of. In fact, we are unaware of any well-posedness results for this equation aside from the 2D Euler case in \cite{BFM2016}, a gap in the literature which should be filled. Accordingly, we will not attempt to rigorously justify equation \eqref{eq:noiseME} in this article. Instead, we will only show in this section how to use our methods to estimate the last term in \eqref{eq:StI} \emph{pathwise} in the noise. Since the sub-Coulombic case is strictly easier than the Coulomb case, we expect that the results of \cite{BFM2016} can be generalized to this case without issue. When combined with the results of the present paper and the reasoning in \cite{Rosenzweig2020spv}, all computations should be justifiable in the sub-Coulombic case in a straightforward manner. 

The assumptions \eqref{ass0} -- \eqref{ass3a'} in the case $0\leq s\leq d-2$ and \eqref{ass4} -- \eqref{ass9'} in the case $d-2<s<d$ carry over. But we also need to supplement \eqref{ass3'} and \eqref{ass9'}, respectively, with the assumptions that the $(\R^d)^{\otimes (4+2k)}$-valued kernel
\begin{equation}
\k(x-y) \coloneqq (x-y)^{\otimes 2} \otimes \nabla^{\otimes (2+2k)}\g(x-y)
\end{equation}
and the $(\R^d)^{\otimes 6}$-valued kernel
\begin{equation}
\K(X-Y) \coloneqq (X-Y)^{\otimes 2} \otimes \nabla^{\otimes 4}\G(X-Y)
\end{equation}
are associated to Calder\'{o}n-Zygmund operators on $\R^{d}$ and $\R^{d+m}$.

\subsection{Second-order commutator estimate}
Following the strategy of \cref{sec:comm}, we prove a second-order version of the commutator estimate \cref{propcommu}. With more work, one can extend the method of proof to arbitrary $k$-th order, but having no need for such generality, we will not do so in this work.

\begin{proposition}
\label{prop:commso}
Let $0\leq s<d$. Let $\g \in C^\infty(\R^d\setminus\{0\})$ such that $\g(x)=\g(-x)$ in $B(0,r_0)$, and
\begin{equation}\label{eq:socommk}
\forall k\geq 1, \qquad |\nabla^{\otimes k}\g(x)| \lesssim |x|^{-k-s} \qquad \forall x\neq 0.
\end{equation}
If $s\leq d-2$, then also assume that
\begin{equation}\label{eq:socommft}
|\hat{\g}(\xi)| \lesssim |\xi|^{s-d} \qquad \forall \xi\neq 0.
\end{equation}
If $s=d-2k$, for some $k\in\N$, also assume that the $(\R^d)^{\otimes (4+2k)}$-valued kernel
\begin{equation}\label{eq:socommCZ}
\k(x,y) \coloneqq (x-y)^{\otimes 2} \otimes \nabla^{\otimes (2+ 2k)}\g(x-y)
\end{equation}
is associated to a Calder\'{o}n-Zygmund operator. Let $v$ be a Lipschitz continuous vector field on $\R^d$. There exists a constant $C>0$, depending only on $s,d$, and the potential $\g$ through \eqref{eq:commgk} -- \eqref{eq:commCZker}, such that for any $f,g\in\Sc(\R^d)$,\footnote{If $s=0$, then it is again implicit that the Fourier transforms of $f,g$ vanish sufficiently rapidly at the origin so that the $\dot{H}^{-d/2}$ norm converges.} we have
\begin{multline}\label{eq:commso}
\left|\int_{(\R^d)^2} (v(x)-v(y))^{\otimes 2} : \nab^{\otimes 2} \g(x-y) f(x) g(y)dxdy\right|\\ 
\leq C\paren*{\|\nabla v\|_{L^\infty} + \|\Dm^{\frac{d-s}{2}}v\|_{L^{\frac{2d}{d-2-s}}} \indic_{s<d-2}}^2  \|f\|_{\dot{H}^{\frac{s-d}{2}}}\|g\|_{\dot{H}^{\frac{s-d}{2}}}.
\end{multline}
Consequently, the integral in the left-hand side of \eqref{eq:commso} extends to a bounded bilinear form $B_v(\cdot,\cdot)$ on $\dot{H}^{\frac{s-d}{2}}(\R^d)$ satisfying the bound \eqref{eq:commso}.
\end{proposition}
\begin{proof}
We follow the outline of the proof of \cref{propcommu}. Since the arguments are very similar, we only sketch the proof, focusing on what is different. Again, we will not explicitly track the dependence of the implicit constants.

\medskip
\noindent
{\bf The super-Coulombic case}. Recycling the kernel notation $\k_v$, define
\begin{equation}
\k_v(x,y) \coloneqq (v(x)-v(y))^{\otimes 2}: \nabla^{\otimes 2}\g(x-y) \qquad \forall x\neq y.
\end{equation}
Note that by assumption \eqref{ass7} and the mean-value theorem, we have the kernel estimate
\begin{equation}
|\k_v(x,y)| \leq C\min\paren*{\frac{\|\nabla v\|_{L^\infty}^2}{|x-y|^s}, \frac{\|v\|_{L^\infty}^2}{|x-y|^{s+2}}}.
\end{equation}
By approximation, we assume that $v \in C^\infty$; $f,g\in\Sc(\R^d)$ with Fourier support away from the origin; and $\g$ has compact support. Writing
\begin{equation}
f = \div \nabla\D^{-1} f \eqqcolon \div f_1
\end{equation}
and similarly for $g$, we find
\begin{equation}
\int_{(\R^d)^2} \k_v(x,y)f(x)g(y)dxdy = \lim_{\ep\rightarrow 0^+} \int_{|x-y|\geq\ep} \k_v(x,y)\div f_1(x)\div g_1(y)dxdy.
\end{equation}
Integrating by parts and using assumption \eqref{eq:socommk} to estimate the boundary terms as in the proof of \cref{propcommu}, we find that
\begin{multline}\nonumber
\int_{|x-y|\geq\ep} \!\!\!\k_v(x,y) f(x)g(y)dxdy = O(\|\nabla^{\otimes 2} v\|_{L^\infty}^2\|fg_1\|_{L^1}\ep^{d-s}) 
+ O(\|v\|_{W^{2,\infty}}^2\|\nabla f_1 g_1\|_{L^1}\ep^{d-s})  
\\
+ O(\|v\|_{W^{1,\infty}}^2 \|\nabla^{\otimes 2}f_1\|_{L^\infty}\|g_1\|_{L^1}\ep^{d-s}) 
 -\frac{1}{2}\int_{|x-y|\geq\ep} \!\!\!\!\!\!\!\!\!\!\!\!\!\!\!\!\!\nabla_x\nabla_y \k_v(x,y) : \paren*{(f_1(x)-f_1(y)) \otimes (g_1(x)-g_1(y))}dxdy.
\end{multline}
The first three terms of the right-hand side vanish as $\ep\rightarrow 0$ since $s<d$. For the last term, direct computation reveals
\begin{equation}
\label{eq:kvsd}
\begin{split}
\nabla_x\nabla_y \k_v(x,y) &= -\nabla^{\otimes 4}\g(x-y) : (v(x)-v(y))^{\otimes 2}\\
&\ph - 2\nabla^{\otimes 3}\g(x-y) : \paren*{(\nabla v(x) + \nabla v(y))\otimes (v(x)-v(y))} \\
&\ph -2\nabla^{\otimes 2}\g(x-y) : (\nabla v(x)\otimes \nabla v(y)).
\end{split}
\end{equation}
So using the triangle inequality and mean-value theorem, we find that
\begin{equation}
\left|\nabla_x\nabla_y \k_v(x,y) \right| \lesssim \frac{\|\nabla v\|_{L^\infty}^2}{|x-y|^{s+2}},
\end{equation}
which implies by Cauchy-Schwarz that
\begin{align}
&\int_{|x-y|\geq\ep} |\nabla_x\nabla_y \k_v(x,y)| |f_1(x)-f_1(y)||g_1(x)-g_1(y)|dxdy \nn\\
&\lesssim \|\nabla v\|_{L^\infty}^2\paren*{\int_{(\R^d)^2} \frac{|f_1(x)-f_1(y)|^2}{|x-y|^{s+2}}dxdy}^{1/2}\paren*{\int_{(\R^d)^2} \frac{|g_1(x)-g_1(y)|^2}{|x-y|^{s+2}}dxdy}^{1/2} \nn\\
&\lesssim\|\nabla v\|_{L^\infty}^2 \|f\|_{\dot{H}^{\frac{s-d}{2}}} \|g\|_{\dot{H}^{\frac{s-d}{2}}}.
\end{align}
In conclusion,
\begin{equation}
\begin{split}
\lim_{\ep\rightarrow 0} \left|\int_{|x-y|\geq\ep} \k_v(x,y) f(x)g(y)dxdy\right| \lesssim \|\nabla v\|_{L^\infty}^2 \|f\|_{\dot{H}^{\frac{s-d}{2}}} \|g\|_{\dot{H}^{\frac{s-d}{2}}}.
\end{split}
\end{equation}

\vspace{0.5cm}
\noindent
{\bf The Coulombic case}. To modify the proof of \cref{propcommu} in the Coulombic case $s=d-2$, we first note that if $|x-y|=\ep$, then writing $x=y+ z$,
\begin{align}
\frac{(x-y)^{i_1}}{|x-y|}\p_{y_{i_2}}\k_v(x,y) &=\frac{z^{i_1}}{|z|} \p_{y_{i_2}}\k_v(y+z,y) \nn\\
&= -\frac{z^{i_1}}{|z|}  \paren*{(z^{i_3}\p_{i_3}v)^{\otimes 2} : \nabla^{\otimes 2}\p_{i_2}\g(z)  + 2(\p_{i_2} v(y) \otimes (z^{i_3}\p_{i_3} v(y))) : \nabla^{\otimes 2}\g(z)} \nn\\
&\ph + O(\|v\|_{W^{2,\infty}}^2\ep^d),
\end{align}
as $\ep\rightarrow 0$, where to obtain the error bound we use assumption \eqref{eq:socommk}. So by Cauchy-Schwarz and assumption \eqref{eq:socommk},
\begin{align}
\left|\int_{\R^d} (f_1(y)\otimes g_1(y)) : \int_{\p B(0,\ep)} \frac{z}{|z|} \otimes \nabla_y\k_v(y+z,y) d\H^{d-1}(z)dy\right| &\lesssim \|\nabla v\|_{L^\infty}^2 \|f_1\|_{L^2} \|g_1\|_{L^2} \nn\\
&\lesssim \|\nabla v\|_{L^\infty}^2 \|f\|_{\dot{H}^{-1}} \|g_1\|_{\dot{H}^{-1}}.
\end{align}
This takes care of the first boundary term in the identity \eqref{eq:Couibp}. The second boundary term vanishes as $\ep\rightarrow 0$ by modifying \eqref{eq:commbt}, \eqref{eq:commbt'} to account for the new definition of $\k_v$.

For the non-boundary term in \eqref{eq:Couibp}, we see, using the identity \eqref{eq:kvsd}, that it is equal to
\begin{equation}\label{eq:sonbd}
\begin{split}
&-\int_{|x-y|\geq\ep}\nabla^{\otimes 2}\p_{i_1i_2}\g(x-y): (v(x)-v(y))^{\otimes 2} f_1^{i_1}(x)g_1^{i_2}(y)dxdy \\
&-2\int_{|x-y|\geq\ep} \nabla^{\otimes 2}\p_{i_2}\g(x-y) : \p_{i_1}v(x) \otimes (v(x)-v(y))f_1^{i_1}(x)g_1^{i_2}(y)dxdy\\
&-2\int_{|x-y|\geq\ep} \nabla^{\otimes 2}\p_{i_1}\g(x-y) : (v(x)-v(y))\otimes  \p_{i_2}v(y) f_1^{i_1}(x)g_1^{i_2}(y)dxdy\\
&-2\int_{|x-y|\geq\ep} \nabla^{\otimes 2}\g(x-y) : (\p_{i_1}v(x) \otimes \p_{i_2}v(y))f_1^{i_1}(x)g_1^{i_2}(y)dxdy.
\end{split}
\end{equation}
By Taylor's theorem,
\begin{align}
&\p_{i_1}v(x) \otimes (v(x)-v(y)):\nabla^{\otimes 2}\p_{i_2}\g(x-y) \nn\\
&= \p_{i_1}v(x)\otimes\paren*{\int_0^1 \nabla v(y+t(x-y))\cdot(x-y)dt} : \nabla^{\otimes 2}\p_{i_2}\g(x-y) \nn\\
&\eqqcolon \p_{i_1}v(x)\otimes \paren*{\int_0^1 \p_{i_3}v(y+t(x-y))dt} : \k_{i_2}^{i_3}(x-y),
\end{align}
where for fixed $i_2,i_3$, the $(\R^d)^{\otimes 2}$-valued kernel $\k_{i_2}^{i_3}(x-y)$ is defined by
\begin{equation}
\k_{i_2}^{i_3}(x-y) \coloneqq (x-y)^{i_3}\nabla^{\otimes 2}\p_{i_2}\g(x-y), \qquad \forall x\neq y.
\end{equation}
By assumption \eqref{eq:socommCZ}, $\k_{i_2}^{i_3}$ is associated to a Calder\'{o}n-Zygmund operator. Hence, the expression
\begin{equation}
\bigg(\int_0^1 \p_{i_3}v(y+t(x-y))dt\bigg)\cdot \k_{i_2}^{i_3}(x-y)
\end{equation}
is the kernel of an $(\R^d)^{\otimes 2}$-valued Calder\'{o}n $d$-commutator, which by the Christ-Journ\'{e} theorem \cite{CJ1987} is bounded on $L^2(\R^d)$. So by Cauchy-Schwarz,
\begin{align}
&\Bigg|\int_{|x-y|\geq\ep} \nabla^{\otimes 2}\p_{i_2}\g(x-y) : \p_{i_1}v(x) \otimes (v(x)-v(y)) f_1^{i_1}(x)g_1^{i_2}(y)dxdy\Bigg|\nn\\
&\lesssim \|\nabla v\|_{L^\infty} \|\p_{i_1}v \otimes f_1^{i_1}\|_{L^2} \|g_1\|_{L^2} \nn\lesssim \|\nabla v\|_{L^\infty}^2 \|f_1\|_{L^2}\|g_1\|_{L^2} \nn
\lesssim \|\nabla v\|_{L^\infty}^2 \|f\|_{\dot{H}^{-1}}\|g\|_{\dot{H}^{-1}}.
\end{align}
By symmetry under swapping the $i_1 \leftrightarrow i_2$, this takes care of the second and third lines of \eqref{eq:sonbd}. For the fourth line, similarly using Taylor's theorem yields
\begin{align}
&\paren*{v(x)-v(y)}^{\otimes2} : \nabla^{\otimes 2}\p_{i_1i_2}\g(x-y) \nn\\
&= \paren*{\int_0^1 \nabla v(y + t(x-y))\cdot (x-y)dt}^{\otimes 2} : \nabla^{\otimes 2}\p_{i_1i_2}\g(x-y) \nn\\
&\eqqcolon \paren*{\int_0^1 \p_{i_3} v(y+t(x-y))dt} \otimes \paren*{\int_0^1 \p_{i_4} v(y+t(x-y))dt}  : \k_{i_1i_2}^{i_3i_4}(x-y), \label{eq:soCald}
\end{align}
where the $(\R^d)^{\otimes 2}$-valued kernel $\k_{i_1i_2}^{i_3i_4}$ is defined by
\begin{equation}
\k_{i_1i_2}^{i_3i_4}(x-y) \coloneqq (x-y)^{\otimes 2,i_3i_4} \nabla^{\otimes 2}\p_{i_1i_2}\g(x-y) \qquad \forall x\neq y.
\end{equation}
By assumption \eqref{eq:socommCZ}, $\k_{i_1i_2}^{i_3i_4}$ is associated to a Calder\'{o}n-Zygmund operator. Hence, the expression \eqref{eq:soCald} is the kernel of a Calder\'{o}n $d$-commutator, and we can apply the Christ-Journ\'{e} theorem \cite{CJ1987} like before to conclude that the modulus of the fourth line of \eqref{eq:sonbd} is bounded by
$\|\nabla v\|_{L^\infty}^2 \|f\|_{\dot{H}^{-1}}\|g\|_{\dot{H}^{-1}}$.

\vspace{0.5cm}

\noindent
{\bf The sub-Coulombic case}. Finally, in the sub-Coulombic case $0\leq s<d-2$, the induction argument proceeds as in the proof of \cref{propcommu}. Instead of identity \eqref{eq:sect}, we now have
\begin{equation}\label{eq:sosect}
\begin{split}
&\int_{(\R^d)^2}\!\!\!\!\!\!\k_v(x,y)f(x)g(y)dxdy= \lim_{\ep\rightarrow0}\Bigg(-\!\!\int_{|x-y|\geq \ep}\!\!\! \!\!\!\! \!\!\!\!\!\!\!\nabla^{\otimes 2}\p_{i_1i_2}\g(x-y) : \paren*{v(x)-v(y)}^{\otimes 2} f_1^{i_1}(x)g_1^{i_2}(y)dxdy \\
&\ph-2\int_{|x-y|\geq \ep} \nabla^{\otimes 2}\p_{i_2}\g(x-y) : \p_{i_1}v(x) \otimes \paren*{v(x)-v(y)}f_1^{i_1}(x)g_1^{i_2}(y)dxdy\\
&\ph-2\int_{|x-y|\geq\ep} \nabla^{\otimes 2}\p_{i_1}\g(x-y) : \paren*{v(x)-v(y)}\otimes \p_{i_2}v(y)  f_1^{i_1}(x)g_1^{i_2}(y)dxdy \\
&\ph- \int_{|x-y|\geq \ep} \nabla^{\otimes 2}\g(x-y) : \paren*{\p_{i_1}v(x)\otimes \p_{i_2}v(y) + \p_{i_2}v(y)\otimes \p_{i_1}v(x)} f_1^{i_1}(x)g_1^{i_2}(y) dxdy \Bigg).
\end{split}
\end{equation}
Fix $i_1,i_2$. For the first term, we can repeat repeat the proof above for the super-Coulombic ($d-4<s<d-2$) and Coulombic ($s=d-4$) cases with $\g$ replaced by $\p_{i_1i_2}\g$. The second and third terms are symmetric, so it suffices to consider the second term. Writing
\begin{align}
\nabla^{\otimes 2}\p_{i_2}\g(x-y) : \p_{i_1}v(x)\otimes \paren*{v(x)-v(y)} = \p_{j_2}\p_{j_1 i_2}\g(x-y) \p_{i_1}v^{j_1}(x) \paren*{v(x)-v(y)}^{j_2} \nn\\
= \nabla\p_{j_1i_2}\g(x-y)\cdot\paren*{v(x)-v(y)}\p_{i_1}v^{j_1}(x)
\end{align}
and therefore
\begin{equation}
\begin{split}
&\lim_{\ep\rightarrow 0} \int_{|x-y|\geq \ep} \nabla^{\otimes 2}\p_{i_2}\g(x-y) : \p_{i_1}v(x) \otimes \paren*{v(x)-v(y)}f_1^{i_1}(x)g_1^{i_2}(y)dxdy \\
&=\lim_{\ep\rightarrow 0} \int_{|x-y|\geq \ep} \nabla\p_{j_1i_2}\g(x-y)\cdot\paren*{v(x)-v(y)} \p_{i_1}v^{j_1}(x)f_1^{i_1}(x)g_1^{i_2}(y)dxdy,
\end{split}
\end{equation}
we see that we can apply \cref{propcommu} to this expression with $\g$ replaced by $\p_{j_1i_2}\g$, $f$ replaced by $\p_{i_1}v^{j_1}f_1^{i_1}$ and $g$ replaced by $g_1^{i_2}$. Thus, its modulus is bounded by
\begin{equation}
\|\nabla v\|_{L^\infty} \|\p_{i_1}v^{j_1}f_1^{i_1}\|_{L^2} \|g_1^{i_2}\|_{L^2} \lesssim \|\nabla v\|_{L^\infty}^2 \|f\|_{\dot{H}^{-1}}\|g\|_{\dot{H}^{-1}}.
\end{equation}
Finally, for the fourth term, $\nabla^{\otimes 2}\g$ is locally integrable, so we can use dominated convergence to pass to the limit $\ep\rightarrow 0$. By symmetry under swapping $x\leftrightarrow y$, it suffices to consider the expression
\begin{equation}
\int_{(\R^d)^2} \nabla^{\otimes 2}\g(x-y) : \paren*{\p_{i_1}v(x)\otimes \p_{i_2}v(y)} f_1^{i_1}(x)g_1^{i_2}(y) dxdy.
\end{equation}
Writing $\mathrm{I} = \Dm^{\frac{d-s}{2}-1}\Dm^{\frac{s-d}{2}+1}$ and using Cauchy-Schwarz, the modulus of the preceding expression is bounded by
\begin{equation}
\|\Dm^{\frac{d-s}{2}-1}(\nabla^{\otimes 2}\g\ast(\p_{i_2}v g_1^{i_2}))\|_{L^2} \|\Dm^{\frac{s-d}{2}+1}( \p_{i_1}v f_1^{i_1})\|_{L^2}.
\end{equation}
By Plancherel's theorem and assumption \eqref{eq:socommft},
\begin{equation}
\|\Dm^{\frac{d-s}{2}-1}(\nabla^{\otimes 2}\g\ast(\p_{i_2}v g_1^{i_2}))\|_{L^2} \lesssim \|\Dm^{\frac{s-d}{2}+1}(\p_{i_2}v g_1^{i_2})\|_{L^2}.
\end{equation}
Using \eqref{eq:dualPlan} and \eqref{eq:dualFL}, we conclude that
\begin{equation}\nonumber
\begin{split}
&\Bigg|\int_{(\R^d)^2} \nabla^{\otimes 2}\g(x-y) : \paren*{\p_{i_1}v(x)\otimes \p_{i_2}v(y)} f_1^{i_1}(x)g_1^{i_2}(y) dxdy\Bigg| \\
&\lesssim \paren*{\|\nabla v\|_{L^\infty}+ \|\Dm^{\frac{d-s}{2}}v\|_{L^{\frac{2d}{d-2-s}}}}^2 \|f\|_{\dot{H}^{\frac{s-d}{2}}}\|g\|_{\dot{H}^{\frac{s-d}{2}}}.
\end{split}
\end{equation}
Just as before, the general case $d-2(k+1)\leq s<d-2k$ follows by induction on $k$, and the existence of the form extension $B_{v}$ by density of $\Sc$ in $\dot{H}^{\frac{s-d}{2}}$.
\end{proof}

\subsection{Renormalization of second-order commutator estimate}
We now use the smearing procedure to obtain a renormalization of \cref{prop:commso}, which is the second-order analogue of \cref{prop:rcomm}.

\begin{proposition}
\label{prop:rcommso}
Let $d\geq 1$, $0\leq s<d$. Let $\ux_N\in (\R^d)^N$ be pairwise distinct, and $\mu \in \P(\R^d)\cap L^\infty(\R^d)$.  If $s=0$, assume that $\int_{\R^d}\log(1+|x|)d\mu(x)<\infty$ and if $s\geq d-1$, assume that \\ $\|\Dm^{s-d}\nabla\mu\|_{L^\infty}<\infty$. Let be $v$ a continuous vector field on $\R^d$ such that
\begin{equation}
\begin{cases}
\|\nabla v\|_{L^\infty} + \|\Dm^{\frac{d-s}{2}}v\|_{L^{\frac{2d}{d-2-s}}} \indic_{s<d-2}< \infty, & {0\leq s\leq d-2} \\
\|v\|_{C^{1,\al}} + \|v\|_{ W^{\frac{d(d+m-s)+2m}{2(d+m)}, \frac{2(d+m)}{d+m-2-s}}} < \infty, & {d-2<s<d},
\end{cases}
\end{equation}
for some $0<\al\leq 1$ and where $m$ is the dimension of the extension space. There exists a constant $C$ depending only $d,s,\al,m$ and on the potential $\g$ through assumptions \eqref{ass1} - \eqref{ass9'}, such that
\begin{multline}
\Bigg| \int_{(\R^d)^2\setminus\triangle} (v(x)-v(y))^{\otimes 2}\cdot\nabla^{\otimes 2}\g(x-y)d\Bigg(\frac{1}{N}\sum_{i=1}^N\d_{x_i}-\mu\Bigg)^{\otimes 2}(x,y)\Bigg| \\
\leq C\Bigg(\|\nabla v\|_{L^\infty}+\|\Dm^{\frac{d-s}{2}}v\|_{L^{\frac{2d}{d-2-s}}}\indic_{s<d-2} 
+\paren*{\|v\|_{C^{1,\al}} + \|v\|_{W^{\frac{d(d+m-s)+2m}{2(d+m)}, \frac{2(d+m)}{d+m-2-s}}}}\indic_{s>d-2}\Bigg)^2\\
\times\Bigg(\Fr_N(\ux_N,\mu) +C(1+\|\mu\|_{L^\infty})\Big(N^{-\frac{\min\{d-s,2\}}{\min\{d,s+2\}(s+1)}} + (N^{-1}\log N) \indic_{s=0} 
+ (N^{-\frac{2}{d(d-1)}}\log N)\indic_{s=d-2}\Big)\Bigg) \\
+ \Bigg(N^{-\frac{\min\{d+1,s+3\}}{\min\{d,s+2\}(s+1)}}\|v\|_{C^{1,\al}}^2 \paren*{\|\Dm^{s+1-d}\mu\|_{L^\infty}\indic_{s<d-1} + \|\Dm^{s-d}\nabla\mu\|_{L^\infty}\indic_{s\geq d-1}} \\
 + \|v\|_{C^{1,\al}}\|\div v\|_{L^\infty}\Bigg(\|\mu\|_{L^\infty}N^{-\frac{(d-s)\min\{d+1,s+3\}}{\min\{d,s+2\}(s+1)}} + N^{-\frac{\min\{d+1,s+3\}}{\min\{d,s+2\}(s+1)}}\Big(\|\mu\|_{L^\infty}\log N + 1\Big)\indic_{s=d-1}\Bigg)\Bigg).
\end{multline}
\end{proposition}
\begin{proof}
We follow the outline of the proof of \cref{prop:rcomm}, focusing on what is different in the second-order setting. Let $\tl{v}$ be the same extension of the vector field $v$ as before. Recycling the kernel notation $\K_v$, we define
\begin{equation}
\label{eq:sodefKv}
\K_v(X,Y) \coloneqq \nabla^{\otimes 2}\G(X-Y) : (\tl{v}(X)-\tl{v}(Y))^{\otimes 2} \qquad \forall X\neq Y\in \R^{d+m}.
\end{equation}
We note that from  assumption \eqref{ass7} that
\be \label{soborneK}
|\K_v(X,Y) |\le  C\|\nab\tl{v}\|_{L^\infty}^2\paren*{|X-Y|^{-s}\indic_{s>0} + \indic_{s=0}}\ee
and 
\be \label{soborneDK} 
|\nab_X \K_v(X,Y)|\le C \|\nab \tilde v\|_{L^\infty}^2 |X-Y|^{-s-1}. \ee
We again let $X_i$ denote $(x_i,0)$.

\medskip

Adding and subtracting the smeared point mass $\d_{X_i}^{(\eta)}$, we obtain the decomposition
\begin{equation}
\int_{(\R^d)^2\setminus\triangle} \k_v(x,y)d\paren*{\frac{1}{N}\sum_{i=1}^N\d_{x_i}-\mu}^{\otimes 2}(x,y) = \Te_1+\Te_2+\Te_3,
\end{equation}
where $\Te_1,\ldots,\Te_3$ are defined as before. We proceed to estimate each of the $\Te_j$ individually.

\noindent{\bf Estimate for $\Te_1$}.
We apply \cref{prop:commso} with
\begin{equation}
f = g=  \frac{1}{N}\sum_{i=1}^N\d_{X_i}^{(\eta)}-\tl{\mu} \in \dot{H}^{\frac{s-d-m}{2}}(\R^{d+m})
\end{equation}
followed by \cref{prop:MElb} or \cref{prop:extMElb} to obtain
\begin{multline}\label{soterm1concl}
|\Te_1|\le C \paren*{\|\nabla\tl{v}\|_{L^\infty} + \|\Dm^{\frac{d+m-s}{2}}\tl{v}\|_{L^{\frac{2(d+m)}{d+m-2-s}}} \indic_{s<d+m-2} }^2  
\\
\times\( \Fr_N(\ux_N, \mu) + \frac{\G_\eta(0)}{N} + C\|\mu\|_{L^\infty}\paren*{\eta^{d-s} + (\eta^d|\log\eta|)\indic_{s=0} + (\eta^2|\log\eta|)\indic_{s=d-2}} +C\eta^2\).
\end{multline}

\noindent{\bf Estimate for $\Te_2$}.
Suppose first that $0\leq s<d-1$. We break $\Te_2$ into  a sum over the index $i$ of terms 
\begin{equation}
\int_{\R^{d+m}}\int_{\R^{d+m}} \paren*{\K_v(X,X_i)-\K_v(X,Y)}d\d_{X_i}^{(\eta)}(Y)d\tl{\mu}(X).
 \end{equation}
Unpacking the definition \eqref{eq:sodefKv} of $\K_v$, we have 
\begin{multline}
\K_v(X,X_i )- \K_v(X,Y) = \paren*{\nabla^{\otimes 2}\G(X-X_i) - \nabla^{\otimes 2}\G(X-Y)} : \paren*{(\tl{v}(X)-\tl{v}(Y))\otimes (\tl{v}(X)-\tl{v}(X_i))} \\
+\nabla^{\otimes 2}\G(X-Y) : \paren*{(\tl{v}(X)-\tl{v}(Y))\otimes (\tl{v}(Y)-\tl{v}(X_i))}\\
 + \nabla^{\otimes 2}\G(X-X_i) : \paren*{(\tl{v}(Y)-\tl{v}(X_i))\otimes (\tl{v}(X)-\tl{v}(X_i))}.
\end{multline}
By considering the cases $|X-X_i|\leq 2\eta$ and $|X-X_i|>2\eta$, it follows from using \eqref{ass7}, the mean-value theorem, and triangle inequality that
\begin{equation}
|\K_v(X,X_i )- \K_v(X,Y)| \leq C\eta\|\nabla\tl{v}\|_{L^\infty}^2\paren*{|X-X_i|^{-s-1} + |Y-X_i|^{-s-1}}
\end{equation}
for all $Y\in \supp(\d_{X_i}^{(\eta)})$. So, we find that
\begin{equation}
\label{eq:soT20fin}
|\Te_2|\indic_{s<d-1} \leq C\eta\|\nab \tilde v\|_{L^\infty}^2 \|\Dm^{s+1-d}\mu\|_{L^\infty}.
\end{equation}

Now if $d-1\leq s<d$, we integrate by parts in $X$ to write
\begin{multline}\label{T2rhs}
\int_{(\R^{d+m})^2} \K_v(X,Y)d\paren*{\d_{X_i}-\d_{X_i}^{(\eta)}}(Y)d\tl\mu(X)\\
= -\int_{(\R^{d+m})^2}\tl\K_v(X,Y)\cdot d\paren*{\d_{X_i}-\d_{X_i}^{(\eta)}}(Y)d\wt{\nabla\mu}(X) \\
-2\int_{(\R^{d+m})^2}\div\tl{v}(X) \nabla\G(X-Y)\cdot\paren*{\tl{v}(X)-\tl{v}(Y)}d\paren*{\d_{X_i}-\d_{X_i}^{(\eta)}}(Y)d\tl\mu(X),
\end{multline}
where we have defined the vector-valued kernel
\begin{equation}
\tl\K_v(X,Y) \coloneqq \nabla\G(X-Y)\cdot\paren*{\tl{v}(X)-\tl{v}(Y)}^{\otimes 2}.
\end{equation}

Proceeding as in the case $0\leq s<d-1$, the first term in the right-hand side of \eqref{T2rhs} is bounded by
\begin{equation}\label{eq:soT21fin}
C\eta\|\nabla\tl{v}\|_{L^\infty}^2\|\Dm^{s-d}\nabla\mu\|_{L^\infty}.
\end{equation}

For the second term, we observe that
\begin{align}
&\int_{\R^{d+m}}\paren*{\tl{v}(X)-\tl{v}(Y)}\cdot\nabla G(X-Y)d(\d_{X_i}-\d_{X_i}^{(\eta)})(Y) \nn\\
&= \int_{\R^{d+m}}\paren*{\paren*{\tl{v}(X)-\tl{v}(X_i)}\cdot\nabla\G(X-X_i) - \paren*{\tl{v}(X)-\tl{v}(Y)}\cdot\nabla G(X-Y)}d\d_{X_i}^{(\eta)}(Y) \nn\\
&=\paren*{\tl{v}(X)-\tl{v}(X_i)}\cdot\paren*{\nabla\G(X-X_i)-\nabla\G_\eta(X-X_i)}\nn\\
&\ph - \int_{\R^{d+m}}\paren*{\tl{v}(X_i)-\tl{v}(Y)}\cdot\nabla\G(X-Y)d\d_{X_i}^{(\eta)}(Y).
\end{align}
If $|X-X_i|\leq 2\eta$, then it follows from the second line using the mean-value theorem and assumption \eqref{ass7} that
\begin{multline}\label{eq:lqeta}
\left|\int_{\R^{d+m}}\paren*{\tl{v}(X)-\tl{v}(Y)}\cdot\nabla G(X-Y)d(\d_{X_i}-\d_{X_i}^{(\eta)})(Y)\right|\\
\leq C\|\nabla\tl{v}\|_{L^\infty} \paren*{|X-X_i|^{-s} + \int_{\R^{d+m}}|X-Y|^{-s}d\d_{X_i}^{(\eta)}(Y)}.
\end{multline}
If $|X-X_i|>2\eta$, then it follows from the third line using the mean-value theorem, assumption \eqref{ass7}, and repeating the proof of \eqref{bdiffgetaext} with $\G$ replaced by $\nabla\G$ that
\begin{multline}\label{eq:gqeta}
\left|\int_{\R^{d+m}}\paren*{\tl{v}(X)-\tl{v}(Y)}\cdot\nabla G(X-Y)d(\d_{X_i}-\d_{X_i}^{(\eta)})(Y)\right| \\
\leq C\|\nabla\tl{v}\|_{L^\infty}\paren*{\eta^2|X-X_i|^{-s-2} + \eta|X-X_i|^{-s-1}},
\end{multline}
where we also use that the reverse triangle inequality and that $\d_{X_i}^{(\eta)}$ is a probability measure. Combining \eqref{eq:lqeta} and \eqref{eq:gqeta}, we obtain
\begin{multline}
\left|\int_{(\R^{d+m})^2}\div\tl{v}(X) \nabla\G(X-Y)\cdot\paren*{\tl{v}(X)-\tl{v}(Y)}d\paren*{\d_{X_i}-\d_{X_i}^{(\eta)}}(Y)d\tl\mu(X)\right| \\
\leq C\|\nabla\tl{v}\|_{L^\infty}\|\div\tl{v}\|_{L^\infty}\Bigg(\int_{|X-X_i|\leq 2\eta}|X-X_i|^{-s}d\tl{\mu}(X) +\int_{|X-X_i|\leq2\eta}\int_{\R^{d+m}}|X-Y|^{-s} d\d_{X_i}^{(\eta)}(Y)d\tl{\mu}(X) \\
+\int_{|X-X_i|>2\eta}\paren*{\eta^2|X-X_i|^{-s-2}+\eta|X-X_i|^{-s-1}}d\tl{\mu}(X)\Bigg).
\end{multline}
Evidently,
\begin{equation}
\int_{|X-X_i|\leq 2\eta} |X-X_i|^{-s}d\tl{\mu}(X) = \int_{|x-x_i|\leq 2\eta} |x-x_i|^{-s}d\mu(x) \leq C\eta^{d-s}\|\mu\|_{L^\infty}.
\end{equation}
Similarly, if $|X-X_i|\leq 2\eta$ and $Y\in\supp(\d_{X_i}^{(\eta)})$, then $|X-Y|\leq 3\eta$. So by direct majorization and Fubini-Tonelli,
\begin{align}
\int_{|X-X_i|\leq2\eta}\int_{\R^{d+m}}|X-Y|^{-s} d\d_{X_i}^{(\eta)}(Y)d\tl{\mu}(X) &\leq \int_{\R^{d+m}}\int_{|X-Y|\leq 3\eta} |X-Y|^{-s}d\tl{\mu}(X)d\d_{X_i}^{(\eta)}(Y) \nn\\
&\leq C\eta^{d-s}\|\mu\|_{L^\infty}.
\end{align}
Lastly, if $s>d-1$, then by dilation invariance of Lebesgue measure,
\begin{equation}
\int_{|X-X_i|>2\eta}\paren*{\eta^2|X-X_i|^{-s-2}+\eta|X-X_i|^{-s-1}}d\tl{\mu}(X) \leq C\eta^{d-s}\|\mu\|_{L^\infty}.
\end{equation}
If $s=d-1$, then we make the modification
\begin{align}
\int_{|x-x_i|>2\eta} |x-x_i|^{-d}d\mu(x) &= \Bigg(\int_{1\geq |x-x_i|>2\eta} |x-x_i|^{-d}d\mu(x) + \int_{|x-x_i|>1} |x-x_i|^{-d}d\mu(x)\Bigg) \nn\\
&\leq C\|\mu\|_{L^\infty}|\log\eta| + 1,
\end{align}
where we use that $\mu$ is a probability density. Thus, we shown that
\begin{multline}
\label{eq:soT22fin}
\left|\int_{(\R^{d+m})^2}\div\tl{v}(X) \nabla\G(X-Y)\cdot\paren*{\tl{v}(X)-\tl{v}(Y)}d\paren*{\d_{X_i}-\d_{X_i}^{(\eta)}}(Y)d\tl\mu(X)\right| \\
\leq C\|\nabla\tl{v}\|_{L^\infty}\|\div\tl{v}\|_{L^\infty}\Bigg(\|\mu\|_{L^\infty}\eta^{d-s} + \eta\Big(\|\mu\|_{L^\infty}|\log\eta| + 1\Big)\indic_{s=d-1}\Bigg).
\end{multline}

Putting together the estimates \eqref{eq:soT21fin} and \eqref{eq:soT22fin} with \eqref{eq:soT20fin}, we conclude that
\begin{multline}
\label{eq:sokpT2fin}
|\Te_2| \leq C\Bigg(\eta\|\nab \tilde v\|_{L^\infty}^2 \|\Dm^{s+1-d}\mu\|_{L^\infty}\indic_{s<d-1}+ \eta\|\nabla\tl{v}\|_{L^\infty}^2\|\Dm^{s-d}\nabla\mu\|_{L^\infty}\indic_{s\geq d-1} \\
 + \|\nabla\tl{v}\|_{L^\infty}\|\div\tl{v}\|_{L^\infty}\Bigg(\|\mu\|_{L^\infty}\eta^{d-s} + \eta\Big(\|\mu\|_{L^\infty}|\log\eta| + 1\Big)\indic_{s=d-1}\Bigg)\Bigg).
\end{multline}

\noindent{\bf Estimate for $\Te_3$}.
This step proceeds as for $\Te_3$ in the proof of \cref{prop:rcomm}, except now using the kernel bounds \eqref{soborneK} and \eqref{soborneDK}. Ultimately, we obtain
\begin{multline}
\label{eq:sokpT3fin}
|\Te_3| \le C  \|\nabla \tl{v}\|_{L^\infty}^2\Bigg(\Fr_N(\ux_N,\mu) + C\Big(\frac{\eta^{-s} + |\log\eta|\indic_{s=0}}{N}  +\frac{\ep^{-s} + |\log\ep|\indic_{s=0}}{N} + \ep^2 \\
+ \|\mu\|_{L^\infty}(\ep^{d-s}+(\ep^d|\log\ep|)\indic_{s=0} + (\ep^2|\log\ep|)\indic_{s=d-2}) + \eta\ep^{-s-1} \Big) \Bigg).
\end{multline}

\noindent
{\bf Conclusion.}
Combining \eqref{term1concl}, \eqref{eq:kpT2fin}, and \eqref{eq:kpT3fin}, we have shown that
\begin{multline}
\label{eq:sorcommrhs}
\left|\int_{(\R^d)^2\setminus\triangle} \k_v(x,y)d\paren*{\frac{1}{N}\sum_{i=1}^N\d_{x_i}-\mu}^{\otimes 2}(x,y)\right|\\
\leq C\paren*{\|\nabla\tl{v}\|_{L^\infty} + \|\Dm^{\frac{d+m-s}{2}}\tl{v}\|_{L^{\frac{2(d+m)}{d+m-2-s}}} \indic_{s<d+m-2} }^2  
\Bigg(\Fr_N(\ux_N, \mu) + \frac{\eta^{-s} + |\log\eta|\indic_{s=0}}{N} \\
+ C\|\mu\|_{L^\infty}\paren*{\eta^{d-s} + (\eta^d|\log\eta|)\indic_{s=0} + (\eta^2|\log\eta|)\indic_{s=d-2}} +C\eta^2 \Bigg) \\
+C\Bigg((\eta\|\nab \tilde v\|_{L^\infty}^2 \|\Dm^{s+1-d}\mu\|_{L^\infty})\indic_{s<d-1}+ (\eta\|\nabla\tl{v}\|_{L^\infty}^2\|\Dm^{s-d}\nabla\mu\|_{L^\infty})\indic_{s\geq d-1} \\
 + \|\nabla\tl{v}\|_{L^\infty}\|\div\tl{v}\|_{L^\infty}\Bigg(\|\mu\|_{L^\infty}\eta^{d-s} + \eta\Big(\|\mu\|_{L^\infty}|\log\eta| + 1\Big)\indic_{s=d-1}\Bigg)\Bigg)
 \end{multline}
 \begin{multline*}
+C  \|\nabla \tl{v}\|_{L^\infty}^2\Bigg(\frac{\ep^{-s} + |\log\ep|\indic_{s=0}}{N} + \ep^2+ \|\mu\|_{L^\infty}(\ep^{d-s}+(\ep^d|\log\ep|)\indic_{s=0} + (\ep^2|\log\ep|)\indic_{s=d-2}) \\
+ \eta\ep^{-s-1} \Bigg).
\end{multline*}
We may now optimize over $\eta$ and $\ep$ exactly as in \eqref{eq:rcommopt}. Using the extension bounds \eqref{eq:tlvbnd} and noting that $\div\tl{v}(x,0) = \div v(x)$, since $\tl{v}$ vanishes in the $(m+1)$-th through $(m+d)$-th components, we arrive at the stated inequality.
\end{proof}


\bibliographystyle{alpha}
\bibliography{ref}
\end{document}